\documentclass[a4paper, 11pt,draft]{article}

\usepackage[T1]{fontenc}
\usepackage[dvipsnames]{xcolor}
\usepackage{mathtools, amssymb, bbm, amsthm, xparse}
\usepackage[draft=false]{hyperref}

\hypersetup
{
	colorlinks = true,
	linkcolor = {red},
	anchorcolor = {black},
	citecolor ={orange},
	filecolor ={cyan},
	menucolor= {red},
	runcolor ={cyan - same as file color},
	urlcolor= {magenta},
}

\numberwithin{equation}{section}

\newtheorem{theo}{Theorem}[section]
\newtheorem{rema}[theo]{Remark}
\newtheorem{lema}[theo]{Lemma}
\newtheorem{defi}[theo]{Definition}

\NewDocumentCommand{\Prob}{O{what} O{\left(} O{\right)}}{\mathbb{P}#2 #1 #3}

\newcommand{\E}{{\mathbb E}}
\newcommand{\1}{{\mathbbm 1}}

\NewDocumentCommand{\pof}{O{} O{}}{\mathbb{P}_{#2}\left(#1\right)}
\newcommand{\iof}[1][A]{\mathbbm{1}_{#1}}
\NewDocumentCommand{\cmv}{O{X} O{Y} O{}}{\E_{\mathbb{#3}}\left(#1\;|\;#2\right)}	
\newcommand{\mb}[1]{\mathbf{#1}}
\newcommand{\mc}[1]{\mathcal{#1}}
\NewDocumentCommand{\norm}{O{p} O{\infty}}{\left\|#1\right\|_{#2}}

\def\a{{\textbf{\cal a}}}
\def\b{{\textbf{\cal b}}}
\def\c{{\textbf{\cal c}}}
\def\d{{\textbf{\cal d}}}

\def\ccX{{\cal X}}
\def\ccB{{\cal B}}
\def\ccF{{\cal F}}
\def\ccN{{\cal N}}
\def\ccG{{\cal G}}
\def\P{\mathbb{P}}
\newcommand{\bor}{\ccB(\ccX)}
\newcommand{\N}{{\mathbb N}}

\newcommand{\R}{{\mathbb R}}
\newcommand{\Z}{{\mathbb Z}}

\usepackage{environ}
\NewEnviron{eqs}{
	\begin{equation}\begin{split}
	\BODY
	\end{split}\end{equation}
}
\NewEnviron{eqsn}{
	\begin{equation*}
	\begin{split}
	\BODY
	\end{split}
	\end{equation*}
}

\begin{document}

\title{General Bernstein-like inequality for additive functionals of Markov chains\footnote{Research
      supported by the National Science Center, Poland, grant no.
      2015/18/E/ST1/00214}}
\author{Micha\l{} Lema\'nczyk\footnote{Institute of Mathematics, University of Warsaw, Banacha 2, 02-097 Warszawa, Poland. miclem1@wp.pl}}
\date{\today}
\maketitle


\begin{abstract}
   Using the renewal approach we prove Bernstein-like inequalities for additive functionals
   of geometrically ergodic Markov chains, thus obtaining counterparts of inequalities for sums
   of independent random variables. The coefficient in the sub-Gaussian part of our estimate is the asymptotic variance of the additive functional, i.e., the variance of the limiting Gaussian variable in the Central Limit Theorem for Markov chains. This refines earlier results by Adamczak and Bednorz, obtained under the additional assumption of strong aperiodicity of the chain.

\end{abstract}

\section{\large{Introduction}}
    Throughout this paper we assume that $\boldsymbol{\Upsilon}=(\Upsilon_n)_{n\in\N}$ is a Markov chain defined on a probability space $(\Omega, \ccF, \P)$, taking values in a measurable (countably generated) space $(\ccX,\ccB)$, with a transition function $P:\ccX\times\ccB\rightarrow [0,1]$. Moreover, we assume that $\boldsymbol{\Upsilon}$ is $\boldsymbol{\psi}$\textbf{-irreducible}, \textbf{aperiodic} and \textbf{admits a unique invariant probability measure} $\pi$. As usual for any initial distribution $\mu$ on $\mc{X}$ we will write $\pof[\boldsymbol{\Upsilon} \in \cdot][\mu]$ for the distribution of the chain with $\Upsilon_0$ distributed according to the measure $\mu$. We will denote by $\delta_x$ the Dirac's mass at $x$ and to shorten the notation we will use $\P_{x}$ instead of $\P_{\delta_x}$.
    
    We say that $\boldsymbol{\Upsilon}$ is \textbf{geometrically ergodic} if there exists a positive number $\rho < 1$ and a real function $G : \ccX \rightarrow \R$ such that for every starting point $x \in \ccX$ and $n \in \N$,
    \begin{equation}\label{E}
       \left\|P^n (x, \cdot) - \pi(\cdot) \right\|_{TV} \leq G(x) \rho^n,
    \end{equation}
    where $\|\cdot\|_{TV}$ denotes the total variation norm of a measure and $P^n(\cdot, \cdot)$ is the $n$-step transition function of the chain. For equivalent conditions, we refer to Chapter 15 of \cite{MT}.

    We will be interested in tail inequalities for sums of random variables of the form
    $$\P_x \left(\left|\sum_{i=0}^{n-1} f(\Upsilon_i)\right|>t\right),$$
    where $f:\ccX \rightarrow \R$ is a measurable real function and $x\in \mc{X}$ is a starting point. 
	Although our main results, stated in Section 4, do not require
    $f$ to be bounded, we give here a version in the bounded case for the sake of simplicity. This
    version will be easier to compare to the Bernstein inequality for bounded random variables stated in Section 2 (cf. Theorem \ref{CBI}). Below for convenience we set $\log(\cdot) = \ln(\cdot \vee e)$, where $\ln(\cdot)$ is the natural logarithm.
    \begin{theo}[Bernstein-like inequality for Markov chains] \label{BBI}
       Let $\boldsymbol{\Upsilon}$ be a geometrically ergodic Markov chain with state space $\ccX$ and let $\pi$ be its unique stationary probability measure.  Moreover, let $f\colon \ccX \to \R$ be a bounded measurable function such that $\E_\pi f=0$. Furthermore, let $x\in\ccX $. Then we can find constants $K,\tau>0$ depending only on $x$ and the transition probability $P(\cdot,\cdot)$ such that for all $t>0$,

       $$\P_x \left(\left|\sum_{i=0}^{n-1} f(\Upsilon_i)\right|>t\right)\leq
          K \exp\left(-\frac{t^2}{32 n \sigma_{Mrv}^2+  \tau  t \| f \|_\infty  \log n}\right),$$

       \noindent where
       \begin{equation}\label{eq:variance}
          \sigma_{Mrv}^2 = \textnormal{Var}_\pi(f(\Upsilon_0))+2\sum_{i=1}^\infty \textnormal{Cov}_\pi(f(\Upsilon_0),f(\Upsilon_i))
       \end{equation}
       denotes the asymptotic variance of the process $\left(f(\Upsilon_i)\right)_i$.
    \end{theo}
    \begin{rema}\em
       We refer to Theorem \ref{BI2} for a more general counterpart of Theorem \ref{BBI} and to Theorem \ref{SBI} for explicit formulas for $K$ and $\tau$.
    \end{rema}

    Let us comment briefly on the method of proof. We rely on the by now classical regeneration technique of Athreya-Ney and Nummelin (see \cite{AN,N,MT}), which allows to split the sum in question into a random number of one-dependent blocks of random lengths. In the context of tail inequalities this approach has been successfully used, e.g., in  \cite{A,AB,Ciolek_Bertail,BC,C,DGM} and provides Bernstein inequalities of optimal type under an additional assumption of strong aperiodicity of the chain (corresponing to $m=1$ in \eqref{MC} below), which ensures that the blocks are independent and allows for a reduction to inequalities for sums of i.i.d. random variables. However, in the general case the implementation of this method available in the literature leads to loss of correlation structure and as a consequence to suboptimal subgaussian coefficient in Bernstein's inequality (in place of $\sigma^2_{Mrv}$). Our main technical contribution is to propose a regeneration-based approach which allows to preserve the correlation structure and recover the correct asymptotic behaviour, corresponding to the CLT for Markov chains.
    \paragraph{}
    The organization of the article is as follows. After a brief discussion of our results (Section 2) we introduce the notation and provide a short description of the regeneration method (Section 3). Next we state our main theorems at their full strength (Section 4). At the end we present their proofs (Section 7). Along the way we develop auxiliary theorems for one-dependent random variables (Section 5) and bounds on number of regenerations (Section 6). Some technical lemmas concerning exponential Orlicz norms are deferred to Appendix.

\section{Discussion of the main result}\label{section discussion of the main result}

    Let us start by recalling the Bernstein inequality in the i.i.d. bounded case.
    \begin{theo}[Classical Bernstein inequality]\label{CBI}
       If $(\xi_i)_i$ is a sequence of i.i.d. centered random variables such that $\sup_i \|\xi_i \|_{\infty}\leq M$ then
       for $\sigma^2 = \E\xi_i^2$ and any $t > 0$
       \begin{equation*}
          \P\left(\sup_{1\le k\le n}\left|\sum_{i=1}^k \xi_i\right| \geq t\right)\leq 2\exp\left(- \frac{t^2}{2n\sigma^2+ \frac{2}{3} M t}\right).
       \end{equation*}
    \end{theo}

    Let us recall that the CLT for Markov chains (see, e.g., \cite{Chen,N,MT}) guarantees that  under assumptions and notation of Theorem \ref{BBI}
    the sums $\frac{1}{\sqrt{n}} \sum_{i=0}^{n-1} f(\Upsilon_i)$ converge in distribution to the normal distribution $\ccN(0,\sigma^2_{Mrv})$. Thus, the inequality obtained in Theorem \ref{BBI} reflects (up to constants) the asymptotic normal behavior of the sums $\frac{1}{\sqrt{n}}\sum f(\Upsilon_i)$ similarly as  the classical Bernstein inequality in the i.i.d. context. Furthermore, the term $\log n$ which appears in our inequality is necessary. More precisely, one can show that if the following inequality holds for all $t > 0$:
    \begin{equation}\label{eq:model-Bernstein}
       \P_x \left(\left|\sum_{i=0}^{n-1} f(\Upsilon_i)\right|>t\right) \le const \cdot\exp\left(-\frac{t^2}{const \cdot n\sigma^2 + const(x) \cdot a_n t \|f\|_\infty}\right)
    \end{equation}
    for some $a_n = o(n)$ and $\sigma \in \R$ ($const$'s stand for some absolute constants whereas $const(x)$ depends only on $x$ and the Markov chain) then
    one must have $\sigma^2 \ge const \cdot \sigma_{Mrv}^2$. Moreover, it is known that for some geometrically ergodic chains $a_n$ must grow at least logarithmically with $n$ (see \cite{A}, Section 3.3).

    \paragraph{}
    Concentration inequalities for Markov chains and processes have been thoroughly studied in the literature, the (non-comprehensive) list of works concerning this topic includes \cite{A,AB,Ciolek_Bertail,BC,  C,DG,DGM,GWu,JSF,Lez,Lezaud, MPR,MPRII,Pau,Sam,Win}. Some results are devoted to concentration for general functions of the chain (they are usually obtained under various Lipschitz or bounded difference type conditions), others specialize to additive functionals, which are the object of study in our case. Tail inequalities for additive functionals are usually counterparts of Hoeffding or Bernstein inequalities. The former ones do not take into account the variance of the additive functional and are expressed in terms of $\|f\|_\infty$ only. They can be often obtained as special cases of concentration inequalities for general function (see, e.g., \cite{DG,Pau,Sam}). Bernstein type estimates of the form \eqref{eq:model-Bernstein} are considered, e.g., in \cite{A,AB,Ciolek_Bertail, BC,C,DGM,GWu,Lez,Lezaud, MPR,MPRII,Pau,Win} and use various variance proxies $\sigma^2$, which do not necessarily coincide with the limiting variance $\sigma_{Mrv}^2$. In the continuous time case, inequalities of Bernstein type for the natural counterpart the additive functional, involving  asymptotic variance have been obtained under certain spectral gap or Lyapunov type conditions in \cite{GWu, Lez}. For discrete time Markov chains, inequalities obtained in \cite{A,AB,BC,C,DGM} by the regeneration method give \eqref{eq:model-Bernstein} (under various types of ergodicity assumptions and with various parameters $a_n$) with $\sigma^2$, which coincides with $\sigma_{Mrv}^2$ only under additional assumption of strong aperiodicity of the chain. On the other hand the articles \cite{MPR,MPRII,Sam,Win} provide more general results, available for non-necessarily Markovian sequences of random variables, satisfying various types of mixing conditions.  The variance proxies $\sigma^2$ that are used in these references are close to the asymptotic variance, however in general do not coincide with it. For instance the inequality obtained in \cite{MPR}, which is valid in particular for geometrically ergodic chains, uses (in our notation) $\sigma^2 = \textnormal{Var}_\pi(f(\Upsilon_0))+2\sum_{i=1}^\infty |\textnormal{Cov}_\pi(f(\Upsilon_0),f(\Upsilon_i))|$. Comparing with \eqref{eq:variance}, one can see that $\sigma_{Mrv}^2 \le \sigma^2$. In fact one can construct examples when the ratio betweeen the two quantities is arbitrarily large or even $\sigma_{Mrv}^2 = 0$ and $\sigma^2 > 0$. The reference \cite{Win} provides an inequality for \emph{uniformly} geometrically ergodic processes, involving a certain implicitly defined variance proxy $\sigma_n^2$, which may be bounded from above by $\sigma^2$ from \cite{MPR} or by $\textnormal{Var}_\pi(f(\Upsilon_0))+C\|f\|_\infty\E_\pi |f(\Upsilon_0)|$, where $C$ is a constant depending on the mixing properties of the process. For a fixed process, in the non-degenerate situation, when the asymptotic variance is non-zero, it can be substituted for $\sigma_n^2$ at the cost of introducing additional multiplicative constants, depending on the chain and the function $f$.

    To the best of our knowledge Theorem \ref{BBI} is therefore the first tail inequality available for general geometrically ergodic Markov chains (not necessarily strongly aperiodic), which (up to universal constants) reflects the correct limiting Gaussian behavior of additive functionals. The problem of obtaining an inequality of this type was posed in \cite{AB}. Let us remark that quantitative investigation of problems related to the Central Limit Theorems for general aperiodic Markov chains seems to be substantially more difficult than for chains which are strongly aperiodic. For instance optimal strong approximation results are still known only in the latter case \cite{MR}.

\section{Notation and basic properties}
    For any $k, l \in \Z$, $k \le l$ we define \textbf{integer intervals} of consecutive integers
    \begin{equation*}
       [k, l] = \left\{k, k + 1, \ldots, l\right\}, \quad [k, l) = \left\{k, k + 1, \ldots, l - 1\right\}, \quad [k, \infty) = \left\{k, k + 1, \ldots \right\}.
    \end{equation*}
    For any process $\mb{X} = \left(X_i\right)_{i\in \N}$ and $S \subset \N$ we put
    \begin{equation*}
       X_{S} = \left(X_i\right)_{i \in S}, \quad \mc{F}^\mb{X} = \left(\mc{F}_i^\mb{X}\right)_{i\in \N}, \quad \mc{F}_i^\mb{X}= \sigma\left(X_{[0, i]}\right).
    \end{equation*}
    Moreover, for $k \in \N$ we define the corresponding \textbf{vectorized process}
    \begin{equation*}
       \mb{X}^{(k)} = \left(X_i^{(k)}\right)_{i \in \N}, \qquad X_i^{(k)} = X_{[i k,(i + 1)k)}.
    \end{equation*}
    \begin{defi}[Stationarity]
       We say that a process $(X_n)_{n\in \N}$ is {\bf{stationary}} if for any $k \in \N$ the shifted process $(X_{n + k})_{n\in \N}$ has the same distribution as $(X_n)_{n\in \N}$.
    \end{defi}
    \begin{defi}[$m$-dependence]
       Fix $m \in \N$. We say that $(X_n)_{n\in \N}$ is \textbf{$\mb{m}$-dependent} if for any $k \in \N$ the process $(X_n)_{n \le k}$ is independent of the process $(X_n)_{n \geq m + 1 + k}$.
    \end{defi}
    \begin{rema}\em
       Let us note that a process $(X_n)_{n\in \N}$ is $0$-dependent iff the variables $(X_n)_{n\in \N}$ are independent.
       Finally let us give a natural example of a $1$-dependent process $(X_n)_{n\in \N}$. Fix an independent process $(\xi_n)_{n\in \N}$ and a Borel, real function $h:\R^2 \rightarrow \R$. Then $(h(\xi_n, \xi_{n+1}))_{n\in \N}$ is $1$-dependent. Such processes are called two-block factors. It is worth noting that there are one-dependent processes which are not two-block factors (see \cite{BGM}).
    \end{rema}
    \begin{rema} \label{every_mth}\em
       Assume that a process $(X_n)_{n\in \N}$ is $m$-dependent. Then for any $n_0\in\N$ the process $(X_{n_0 + k(m+1)})_{k\in\N}$ is independent. Moreover, if the process  $(X_n)_{n\in \N}$ is stationary then for any $n_0\in\N$, $(X_{n_0 + k(m+1)})_{k\in\N}$ is a collection of i.i.d. random variables.
    \end{rema}
    \subsection{Split chain}\label{section split chain}
        As already mentioned in the introduction our proofs will be based on the regeneration technique which was invented independently by Nummelin and Athreya-Ney (see \cite{AN} and \cite{N}) and was popularized by Meyn and Tweedie \cite{MT}. We will introduce the split chain and then regeneration times of the split chain. The construction of the split chain is well-known and as references we recommend \cite{MT} (Chaps. 5,17) and \cite{N}. We briefly recall this technique bellow. Let us stress that although this construction is based on the one presented in \cite{MT} our notation is slightly different. Firstly let us recall the minorization condition for Markov chains which plays a main role in the splitting technique.

        \begin{defi}\em
           We say that a Markov chain $\boldsymbol{\Upsilon}$ satisfies \textbf{the Minorization Condition} if there exists a set $C\in \bor$ (called a small set),
           a probability measure $\nu$ on $\mc{X}$ (a small measure), a constant $\delta>0$ and a positive integer $m\in\N$ such that  $\pi(C) > 0$ and
           \begin{equation}\label{MC}
              P^m(x,B)\geq \delta \nu(B)
           \end{equation}
           holds for all $x\in C$ and $B\in\bor$.
        \end{defi}
        \begin{rema}\em
           One can assume that $\nu(C)=1$ (possibly at the cost of increasing $m$).
        \end{rema}
        \begin{rema}\em
           One can check that under assumptions of our theorem the Minorization Condition \eqref{MC} holds for some $C$, $\nu$, $\delta$ and $m$. We refer to \cite{MT}, Section 5.2 for the proof of this fact.
        \end{rema}
        Fix $C$, $m$, $\nu$ and $\delta>0$ as in \eqref{MC}. The minorization condition allows us to redefine the chain $\boldsymbol{\Upsilon}$ together with an auxiliary regeneration structure. More precisely we start with a splitting of the space $\ccX$ into two identical copies on level $0$ and $1$ namely we consider $\overline{\ccX} =\ccX\times \{0,1\}$.
        Now we split $\boldsymbol{\Upsilon}$ in the following way. We consider a process $\boldsymbol{\Phi} =(\boldsymbol{\Upsilon}, \boldsymbol{\Lambda})=(\Upsilon_i,\Lambda_i)_{i\geq 0}$  (usually called the split chain) defined on $\overline{\ccX}$ (we slightly abuse the notation by denoting the first coordinate of the split chain with the same letter as for the initial Markov chain, but it will turn out that the first coordinate of the split chain has the same distribution as the starting Markov chain so this notation is justified). The random variables $\Lambda_k$ take values in $\{0,1\}$ (they indicate the level on which  $\Phi_k$ is). For a fixed $x\in C$ let
        \begin{equation}\label{Radon Derivative}
           r(x,y) = \frac{\delta\nu(dy)}{P^m(x,dy)}
        \end{equation}
        and note that the above Radom-Nikodym derivative is well-defined thanks to \eqref{MC}. Moreover, $r(x,y) \leq 1$. Now for any $A_1,\ldots,A_m\in\bor$, $k\in\N$ and $i \in \{0, 1\}$ set
        \begin{equation}\label{Markov-like1}
           \begin{split}
              &\P\left(\Lambda_{km} = i, \Upsilon_{[km+1, (k + 1)m]}\in A_1\times \cdots \times A_m \;|\;\ccF_{km}^{\boldsymbol{\Upsilon}}, \ccF_{km-m}^{\boldsymbol{\Lambda}} ,\Upsilon_{km}=x\right)\\
              &= \P\left(\Lambda_{0} = i, \Upsilon_{[1, m]}\in A_1\times \cdots \times  A_m \;|\; \Upsilon_{0}=x\right)\\
              &=\int_{A_1} \cdots \int_{A_m} r(x,x_m, i)
              P(x_{m-1},dx_m)P(x_{m-2},dx_{m-1})\ldots P(x,dx_{1}),
           \end{split}
        \end{equation}
        where
        \begin{equation}\label{definition of ri}
           r(x, y, i) =
           \begin{dcases}
              \1_{x\in C}\; r(x,x_m), \quad &\textnormal{if } i = 1,  \\
              1-\1_{x\in C}\; r(x,x_m), \quad &\textnormal{if } i = 0. \\
           \end{dcases}
        \end{equation}

        \noindent Moreover, for any $k, i \in \N$ such that $km < i < (k + 1)m$ we set
        \begin{equation}\label{MarkovLikeY}
           \Lambda_i = \Lambda_{km}.
        \end{equation}
        \begin{rema}[Initial distribution for the split chain]\em
           In order to be able to set initial distribution for the split chain for arbitrary probability measure $\mu$ on $\ccX$ we define \textbf{the split measure $\mu^*$} on $\overline{\ccX}$ by:
           \begin{equation}
              \mu^*(A\times \{i\}) =
              \begin{dcases}
                 (1-\delta)\mu(C\cap A)+\mu(A\cap C^c), \quad &\textnormal{if } i = 0, \\
                 \delta\mu(C\cap A), \quad &\textnormal{if } i = 1  .                  \\
              \end{dcases}
           \end{equation}
           Such definition ensures that  $(\Upsilon_0, \Lambda_0) \sim \mu^*$ as soon as $\Upsilon_0 \sim \mu$. For convenience sake, for any $x \in \ccX$, we will write
           \begin{equation}
              \P_{x^*}(\cdot) = \P_{\delta_{x}^*}(\cdot).
           \end{equation}
        \end{rema}

        \begin{rema}[Markov-like properties of the split chain]\label{MarkovPropertySplitChain}\em
           In order to give some intuition behind the definition of the split chain note that the distribution of the first coordinate of the split chain $\boldsymbol{\Phi}$ with initial distribution $\mu^{*}$ coincides with that of the original Markov chain $\boldsymbol{\Upsilon}$ which starts from $\mu$. From now on $\boldsymbol{\Upsilon}$ always corresponds to this first coordinate of the split chain. One can easily generalize \eqref{Markov-like1} to show the following \textbf{Markov-like property of the split chain}: for any $k\in\N$ and product measurable bounded function $F$ we have
           \begin{equation}\label{Markov like property of the split chain}
              \E \left(F\left(\Upsilon_{[km+1, \infty)}, \Lambda_{[km, \infty)}\right)\;|\;\ccF_{km}^{\boldsymbol{\Upsilon}},\ccF_{km-m}^{\boldsymbol{\Lambda}} \right) = \E\left( F\left(\Upsilon_{[km+1, \infty)}, \Lambda_{[km, \infty)}\right)\;|\;\Upsilon_{km} \right).
           \end{equation}
           This, in turn, leads to the fact that the vectorized split chain $\boldsymbol{\Phi}^{(m)}$ \textbf{is a Markov chain}. Even more, for any product measurable bounded function $F$ and $k\in\N$ we have
           \begin{equation*}
              \E\left(F\left(\Phi^{(m)}_{[k, \infty)} \right)|\;\Phi^{(m)}_{[0, k)}\right) = \E\left(F\left(\Phi^{(m)}_{[k, \infty)} \right)|\; \Phi^{(m)}_{k - 1}\right) = \E\left(F\left(\Phi^{(m)}_{[k, \infty)} \right)|\; \Upsilon_{mk-m}, \Upsilon_{mk - 1},  \Lambda_{mk - m}\right).
           \end{equation*}
        \end{rema}
        Now we can introduce the aforementioned \textbf{regeneration structure for} $\boldsymbol{\Phi}$. Firstly we define certain stopping times. For convenience we put $\tau_{-1} = -m$ and then for $i \ge 0$ we define $\tau_i$ to be the $i$'th time when the second coordinate (level coordinate) hits $1$, namely
        \begin{equation}\label{Regeneration times definition}
           \tau_i = \min\{k> \tau_{i-1} \;\big|\; \Lambda_{k}=1, \; m|k\}.
        \end{equation}
        Now we are ready to introduce \textbf{the random blocks and the random block process}
        \begin{equation}
           \Xi_i=\Upsilon_{[\tau_{i-1}+m, \tau_i + m)}, \quad \boldsymbol{\Xi} = \left(\Xi_{i}\right)_{i \ge 0},
        \end{equation}
        where we consider $ \Xi_i$ as a random variable with values in the disjoint union $\bigsqcup_{j \ge 0}\mc{X}^j$. For clarity of this presentation, here and later on, we omit the measurability details.
        \begin{rema}\label{rem:Rndom-blocks}\em
           Let us now briefly discuss the behaviour of these random blocks. Firstly, by the strong Markov property of the vectorized split chain it is not hard to see that $\boldsymbol{\Xi}$ \textbf{is a Markov chain}. On a closer look one can see that for any product measurable function $F$
           \begin{equation}
              \cmv[F\left(\Xi_{[i, \infty)}\right)][\Xi_{[0, i)}] = \cmv[F\left(\Xi_{[i, \infty)}\right)][\Xi_{i - 1}] = \cmv[F\left(\Xi_{[i, \infty)}\right)][\textnormal{pr}_m\left(\Xi_{i - 1}\right)],
           \end{equation}
           where $\textnormal{pr}_m:\bigsqcup_{j \ge m}\mc{X}^j\to \mc{X}^m$ is a projection on $m$-last coordinates,
           \begin{equation}
              \textnormal{pr}_m\left(x_0, \ldots, x_j\right) = \left(x_{j - m + 1}, \ldots, x_j\right).
           \end{equation}
           Apart from being Markovian the sequence $(\Xi_i)_{i \geq 0}$ is $1$-dependent whereas $(\Xi_i)_{i \geq 1}$ is stationary (see \cite{Chen}, Corollary 2.4). The stationarity follows from the fact that for $m|k$ we have
           \begin{equation}
              \mc{L}\left(\Upsilon_{k + m}\;|\; \Lambda_{k} = 1\right) = \nu,
           \end{equation}
           that is every time $k$ (which is a multiple of $m$) the split chain is on level $1$ (note that this implies $\Upsilon_k \in C$) the split chain \textbf{regenerates} and starts anew from $\nu$.
           Furthermore, the lengths of $\Xi_i$,
           \begin{equation}\label{len_block_iid}
              \left|\Xi_i\right| =  \tau_{i} - \tau_{i - 1} \\
           \end{equation}
           are \textbf{independent} random variables for $i \ge 0$ and form a stationary process for $i\ge 1$. Let us add that if $m = 1$, one can show that $\Xi_i$'s are independent. This fact makes a crucial difference between strongly aperiodic and not strongly aperiodic Markov chains (see \cite[Section 6]{BLL}).
        \end{rema}
        At last let us introduce \textbf{the excursions and the excursion process}
        \begin{equation}\label{Excusions Definition}
           \chi_i = \chi_i\left(f\right) = \sum_{j = \tau_i + m}^{\tau_{i + 1}+ m - 1} f(\Upsilon_j), \qquad \boldsymbol{\chi} = \left(\chi_i\right)_{i \ge 0},
        \end{equation}
        which will play a crucial role in our future considerations.
        By properties of the random blocks one concludes that $\boldsymbol{\chi}$ is $1$-dependent and satisfies
        \begin{equation}\label{siconditional}
           \E\left(\chi_i \;|\; \Xi_{[0, i]}\right) = \E\left(\chi_i \;|\; \Xi_{i}\right).
        \end{equation}
        Moreover,  $\left(\chi_i\right)_{i \ge 1}$ is stationary. Due to the Pitman occupation measure formula (see, \cite{MT}, Theorem 17.3.1, page 428) which says that for any measurable real function $G$,
        \begin{equation}\label{POMF}
           \E_\nu \sum_{i=0}^{\tau_0/m} G(\Upsilon_{mi},\Lambda_{mi})=\delta^{-1}\pi(C)^{-1}\E_\pi G(\Upsilon_0,\Lambda_0),
        \end{equation}
        and observation that
        $\P_\mu$-distribution of excursion $\chi_i(f)$ ($i \ge 1$) is equal to the $\P_\nu$ -distribution of
        $\chi_0$, we get that for any initial distribution $\mu$ and any $i \ge 1$,
        \begin{equation}\label{mean_of_excursion}
           \E_\mu \chi_i = \E_\nu \chi_{0}= \delta^{-1}\pi(C)^{-1}m\int f d\pi.
        \end{equation}
        As a consequence, $\E_\pi f(\Upsilon_i) = 0$ implies that for every $i \ge 1$, $\E_\mu \chi_i(f) = 0$.
        Now we are ready to decompose our sums into random blocks. If $m | n$ then
        \begin{equation} \label{block_decomposition}
           \sum_{i=0}^{n-1}f(\Upsilon_i) = \left(\sum_{i=0}^{\tau_0/m} \Theta_i \1_{N>0}  + \1_{N=0} \sum_{i=0}^{n/m - 1} \Theta_i \right)+\left(\sum_{i=1}^{N} \chi_{i-1}(f)\right) - \left( \1_{N>0} \sum_{k=n}^{\tau_N+ m - 1} f(\Upsilon_k)\right),
        \end{equation}
        where
        \begin{equation} \label{eq:definition-of-Theta-N}
           \Theta_k = \Theta_k(f) = \sum_{i=0}^{m-1} f(\Upsilon_{km+i}), \qquad N = \inf\{i\geq 0 \;|\; \tau_i + m - 1 \geq n-1\}.
        \end{equation}
        This decomposition will be of utmost importance in our proof.
    \subsection{Asymptotic variances}
        During the upcoming proofs we will meet two types of asymptotic variances: $\sigma_{Mrv}^2$ associated with the process $\left(f(\Upsilon_i)\right)_{i \ge 0}$ and $\sigma_{
              \infty}^2$ associated with $\boldsymbol{\chi}$. The first one, defined as
        \begin{equation}\label{eq:definition-of-sigma-mrv}
           \sigma_{Mrv}^2 = \lim\limits_{n\to\infty}\frac{1}{n}\textnormal{Var}\left(f(\Upsilon_0) + \cdots + f(\Upsilon_{n - 1}) \right)=   \textnormal{Var}_\pi (f(\Upsilon_0)) + 2 \sum_{i\ge 1} \textnormal{Cov}_\pi(f(\Upsilon_i),f(\Upsilon_0))
        \end{equation}
        is exactly the variance of the limiting normal distribution of the sequence $\frac{1}{\sqrt{n}} \sum_{i=1}^n f(\Upsilon_i)$. The second one,
        \begin{equation*}
           \sigma_\infty^2 = \lim\limits_{n\to\infty}\frac{1}{n}\textnormal{Var}\left(\chi_{1} + \cdots + \chi_n \right) = \E \chi_1^2 + 2\E \chi_1 \chi_2
        \end{equation*}
        is the variance of the limiting normal distribution of the sequence $\frac{1}{\sqrt{n}} \sum_{i=1}^n \chi_i$.
        Both asymptotic variances are very closely linked via the formula
        \begin{equation} \label{asymp_variances_formula}
           \sigma_\infty^2 = \sigma_{Mrv}^2 \E (\tau_1-\tau_0) = \sigma_{Mrv}^2 m \delta^{-1} \pi(C)^{-1}.
        \end{equation}
        For the proof of this formula we refer to \cite{MT} (see $(17.32)$, page 434).

\section{Main results}
    In order to state our results in the general form we need to recall the definition of the exponential Orlicz norm. For any random variable $X$ and $\alpha>0$ we define
    \begin{equation}\label{definition of the exponential Orlicz norm}
       \| X \|_{\psi_\alpha}=\inf\Big\{c>0 \;|\; \E \exp\left(\frac{|X|^\alpha}{c^\alpha}\right)\leq 2\Big\}.
    \end{equation}
    If $\alpha < 1$ then $\|\cdot\|_{\psi_\alpha}$ is just a quasi-norm (for basic properties of these quasi-norms we refer to Appendix \ref{appendix exponential norm}).
    In what follows we will deal with various underlying measures on the state space $\overline{\ccX}$. In order to stress the dependence of the Orlicz norm on the initial distribution $\mu$ of the chain $\boldsymbol{\Phi}$ we will sometimes write $\|\cdot\|_{\psi_\alpha,\mu}$ instead of  $\|\cdot\|_{\psi_\alpha}$.

    Before we formulate our main result let us introduce and explain the role of the following parameters:
    \begin{equation} \label{eq1}
       \a = \left \|\sum_{k=0}^{\tau_0 /m}  \left|\Theta_k\right|\right\|_{\psi_\alpha, \P_{x^*}},  \quad \b = \left\|\sum_{k=0}^{\tau_0 /m}  \left|\Theta_k\right| \right\|_{\psi_\alpha,\P_{\pi^*}},\quad
       \c = \left \|\chi_i(f)\right\|_{\psi_\alpha},\quad  \d  = \left\|\tau_1-\tau_0\right\|_{\psi_1},
    \end{equation}
    where $\Theta_k = \sum_{i=0}^{m-1} f(\Upsilon_{km+i})$ (cf. \eqref{block_decomposition}).
    \noindent The parameter $\a$ (resp. $\b$) will allow us to estimate the first (third) term on the right-hand side of \eqref{block_decomposition}, whereas the parameters $\c$ and $\d$ will be used to control the middle term. We note that  $\d$ quantifies geometric ergodicity of $\boldsymbol{\Upsilon}$ and is finite as soon as $\boldsymbol{\Upsilon}$ is geometrically ergodic. Let us mention that all these parameters can be bounded for example by means of drift conditions widely used in the theory of Markov chains (see Remark \ref{drift-like-conditions}).
    Finally let us remind that $\sigma_{Mrv}^2  = \textnormal{Var}_\pi (f(\Upsilon_0)) +2 \sum_{i = 1}^\infty \textnormal{Cov}_\pi(f(\Upsilon_0),f(\Upsilon_i))$ denotes the asymptotic variance of normalized partial sums of the process $\left(f(\Upsilon_i)\right)_i$.

    We are now ready to formulate the first of our main results (recall the definitions of the small set $C$ and the minorization condition \eqref{MC}).

    \begin{theo} \label{BI}
       Let $\boldsymbol{\Upsilon}$ be a geometrically ergodic Markov chain and $\pi$ be its unique stationary probability measure.  Let $f\colon \ccX \to \R$ be a measurable function such that $\E_\pi f=0$  and let $\alpha\in (0,1]$. Moreover, assume for simplicity that $m|n$. Then for all $x\in\ccX$ and $t>0$,
       \begin{equation}\label{eq:first-theorem}
          \begin{split}
             &\P_x \left(\left|\sum_{i=0}^{n-1} f(\Upsilon_i)\right|>t\right)  \le 2\exp\left(-\frac{t^\alpha}{(23\a)^\alpha}\right) + 2\left[\delta\pi(C)\right]^{-1} \exp\left(-\frac{t^\alpha}{(23\b)^\alpha}\right)\\
             &\;\; +	6\exp(8)\exp\left(-\frac{t^\alpha}{\frac{16}{\alpha} (27 \c)^\alpha}\right)
             + 6\exp\left(-\frac{t^2}{30 n\sigma_{Mrv}^2+ 8 t M}\right) +\exp(1)\exp\left(-\frac{ n m}{67 \delta \pi(C)\d^2 }\right),
          \end{split}
       \end{equation}
       where $\sigma^2_{Mrv}$ denotes the asymptotic variance for the process $(f(\Upsilon_i))_i$ given by  \eqref{eq:definition-of-sigma-mrv}, the parameters $\a,\b,\c,\d$ are defined by \eqref{eq1} and  $M=\c(24\alpha^{-3} \log{n})^\frac{1}{ \alpha}$.
    \end{theo}
    \begin{rema}\em\label{drift-like-conditions}
       For the conditions under which $\a,\b,\c$ are finite we refer to \cite{AB}, where the authors give bounds on $\a,\b,\c$ under classical drift conditions. If $f$ is bounded then one easily shows that
       \begin{equation}\label{c bound}
          \max\left(\a, \b\right) \le 2 D \|f\|_\infty,\qquad \c \leq  D\|f\|_{\infty},
       \end{equation}
       where $D = \max\left(\d, \|\tau_0\|_{\psi_1, \;\P_{x^*}}, \|\tau_0\|_{\psi_1, \;\P_{\pi^*}}\right)$.
       For computable bounds on $D$ we refer to \cite{BAX}.
    \end{rema}

    \noindent Let us note that in Theorem \ref{BI} the right-hand side of the inequality does not converge to $0$ when $t$ tends to infinity (one of the terms depends on $n$ but not on $t$). Usually in applications $t$ is of order at most $n$ and the other terms dominate on the right-hand side of the inequality, so this does not pose a problem. Nevertheless one can obtain another version of Theorem \ref{BI}, namely

    \begin{theo} \label{BI2}
      Under the assumptions and notation of Theorem \ref{BI} we have
       \begin{equation}
          \begin{split}
             &\P_x \left(\left|\sum_{i=0}^{n-1} f(\Upsilon_i)\right|>t\right)  \le 2\exp\left(-\frac{t^\alpha}{(54\a)^\alpha}\right) + 2\left[\delta\pi(C)\right]^{-1} \exp\left(-\frac{t^\alpha}{(54\b)^\alpha}\right)\\
             &\;\; +		4\exp(8)\exp\left(-\frac{t^\alpha}{\frac{16}{\alpha}(27\c)^\alpha}\right) + 6\exp\left(-\frac{t^2}{37(1+p) n\sigma_{Mrv}^2+ 18 M \d \sqrt{K_p}t}\right), \\
          \end{split}
       \end{equation}
       where $K_p = L_p + 16/L_p$  and $L_p =  \frac{16}{p} + 20$.
    \end{theo}
    
    It is well-known that for geometrically ergodic chains $\|\tau_0\|_{\psi_1, \;\P_{x^*}}$, $\|\tau_0\|_{\psi_1, \;\P_{\pi^*}}$, $\|\tau_1 - \tau_0\|_{\psi_1} < \infty$ (see \cite{BAX} for constructive estimates). Therefore \eqref{c bound} and  Theorem \ref{BI} lead to

    \begin{theo} \label{SBI}
       Let $\boldsymbol{\Upsilon}$ be a geometrically ergodic Markov chain and $\pi$ be its unique stationary, probability measure.  Let $f\colon \ccX \to \R$ be a bounded, measurable function such that $\E_\pi f=0$. Fix $x\in\ccX$. Moreover assume that $\|\tau_0\|_{\psi_1,\delta_{x}^*}$, $\|\tau_0\|_{\psi_1,\pi^*}$, $\|\tau_1-\tau_0\|_{\psi_1} \leq D$. Then for all $t>0$,

       \begin{equation}\label{eq:third theorem}
          \begin{split}
             \P_x \left(\left|\sum_{i=0}^{n-1} f(\Upsilon_i)\right|>t\right) &\le K \exp\left(-\frac{t^2}{32 n \sigma_{Mrv}^2+  433 t \delta \pi(C)\|f\|_\infty D^2 \log n} \right),
          \end{split}
       \end{equation}

       \noindent where $\sigma^2_{Mrv}$ is the asymptotic variance of $(f(\Upsilon_i))_i$ and $K = \exp(10)+2\delta^{-1}\pi(C)^{-1}$.
    \end{theo}

    \begin{rema}\em
       Theorem \ref{SBI} implies our main Theorem \ref{BBI} from Introduction with constants $K =  \left(\exp(10) + 2\delta^{-1}\pi(C)^{-1}\right)$ and $\tau =433 \delta \pi(C) D^2 $.
    \end{rema}
\section{Bernstein inequalities for one-dependent sequences}
    In this section we will show two versions (for suprema and randomly stopped sums) of Bernstein inequality for one-dependent random variables. They will be later used in the proofs of our main theorems. In what follows for a one-dependent sequence of random variables $(X_i)_{i\ge 0}$, $\sigma_\infty^2$ denotes the asymptotic variance of normalized partial sums, i.e.,
    \begin{displaymath}
       \sigma_\infty^2 = \E X_1^2 + 2\E X_1X_2.
    \end{displaymath}

    \begin{lema}[Bernstein inequality for suprema of partial sums] \label{1DUBI}
       Let $(X_i)_{i\geq 0}$ be a $1$-dependent sequence of centered random variables such that $\E\exp(c^{-\alpha}|X_i|^\alpha)\leq 2$ for some $\alpha \in (0,1]$ and $c>0$. Assume that there exists a filtration $\left(\ccF_i\right)_{i\ge 0}$ such that for $Z_i = X_i + \E\left(X_{i+1}|\ccF_i\right) -\E\left(X_{i}|\ccF_{i-1}\right)$ we have the following:
       \begin{enumerate}
          \item[0)] $X_i$ is $\ccF_i$ measurable,
          \item[1)] $(Z_i)_{i\ge 1}$ is stationary,
          \item[2)] $(Z_i)_{i\ge 1}$ is $m$-dependent with $m=1$ or $m=2$,
          \item[3)] $\left(\E\left(X_{i}|\ccF_{i-1}\right)\right)_{i\ge 1}$ is stationary,
          \item[4)] $ \E(X_i|\ccF_{i-1})$ is independent of $X_{i+1}$ for any $i \ge 1$.
       \end{enumerate}
       Then
       \begin{equation}\label{asymptotic_variance_formula}
          \E Z_i^2 = \sigma_\infty^2, \quad \|Z_i\|_{\psi_\alpha} \le  c (8/\alpha)^\frac{1}{\alpha}.
       \end{equation}
       Moreover, for any $t > 0$ and $n\in\N$,
       \begin{equation}\label{aux}
          \P\left(\sup_{1 \le k \le n} \left|\sum_{i=1}^k X_i\right| > t \right)  \le K_m\exp\left(-\frac{t^\alpha}{u_m c^\alpha}\right)
          \;\; + L_m\exp\left(-\frac{t^2}{v_{n, m} \sigma_\infty^2+ w_{n, m} t}\right)
       \end{equation}
       where $u_m=\frac{16\cdot 8^\alpha(m+1)^\alpha}{\alpha}$, $v_{n, m}=5(m+1)(n  + m + 1)$, $w_{n, m}=2(m+1)(24\alpha^{-3} \log{n})^\frac{1}{ \alpha}c$, $K_m =  2(m + 1)\exp(8)$ and $L_m = 2(m  + 1)$.
    \end{lema}
    \begin{proof}
       Firstly we will show that if $X_i$'s are centered, independent random variables with common variance $\sigma_\infty^2$ and $\E\exp(c^{-\alpha}|X_i|^\alpha)\leq 2$, then \eqref{aux} holds with $u_0 = 2\cdot 6^\alpha$, $v_{n, 0}=\frac{72}{25}n$, $w_{n, 0}=\frac{8}{5}c \left(3\alpha^{-2}\log n\right)^{\frac{1}{\alpha}}$, $K_0 =  \exp(8)$ and $L_0 = 2$ (allowing for a slight abuse of precision we consider this the $m=0$ case of the lemma).
       Indeed, by Lemma $4.1$ in \cite{AB} for $\lambda = (2^{1/ \alpha} c)^{-1}$,
       \begin{equation}\label{aux2}
          \E\exp\left(\lambda^\alpha\sum_{i=0}^{n-1} \left(|U_i|^\alpha+(\E|U_i|)^\alpha \right)\right)\leq \exp(8),
       \end{equation}
       where $U_i = X_i \iof[\left|X_i\right| > M_0]$ stands for the "unbounded" part of $X_i$ and $M_0 = c \left(3\alpha^{-2}\log n\right)^{\frac{1}{\alpha}}$. Define the "bounded" part of $X_i$, $B_i = X_i \1_{|X_i| \le M_0}$ and notice that $X_i = \overline{B_i} + \overline{U_i}$, where $\overline{B_i} = B_i - \E B_i$ and $\overline{U_i} = U_i - \E U_i$. Using the union bound we get for $p = 1/6$
       \begin{equation*}
          \P\left(\sup_{1\le k \le n}\left|\sum_{i=1}^k X_i\right| > t \right)\le \P\left(\sup_{1\le k \le n}\left|\sum_{i=1}^k \overline{U_i} \right| > tp \right) + \P\left(\sup_{1\le k \le n}\left|\sum_{i=1}^k \overline{B_i} \right| > t(1-p) \right).
       \end{equation*}

       Consider first the unbounded part. Using the
       subadditivity of $x \to x^\alpha$, Markov's inequality and then \eqref{aux2} we get
       \begin{equation*}
          \begin{split}
             &\P\left(\sup_{1\le k \le n}\left|\sum_{i=1}^k \overline{U_i} \right| > tp \right)  \le  \P\left(\exp\left(\lambda^\alpha \sum_{i=1}^n |\overline{U_i}|^\alpha \right) > \exp\left(\lambda p t\right)^\alpha \right)\\
             & \qquad \qquad \le \exp\left(8\right) \exp\left(-\frac{t^\alpha p^\alpha}{2 c^\alpha}\right)  = \exp\left(8\right) \exp\left(-\frac{t^\alpha }{2(6c)^\alpha}\right).
          \end{split}
       \end{equation*}
       As for the "bounded" part, notice that $\E\overline{B_i}^2 \le  \E B_i^2 \le \E X_i^2 = \sigma_\infty^2$. Therefore using the classical Bernstein inequality we get
       \begin{equation*}
          \P\left(\sup_{1\le k \le n}\left|\sum_{i=1}^k \overline{B_i} \right| > t(1-p) \right)  \le 2\exp\left(-\frac{t^2(1-p)^2}{2n\sigma^2_\infty + \frac{4}{3}t(1-p)M_0}\right).
       \end{equation*}
       Combining the three last estimates and substituting $p = 1/6$, allows to finish the proof for independent random variables.

       We will now use the independent case to prove the tail estimate \eqref{aux}, assuming \eqref{asymptotic_variance_formula}, the proof of which we postpone.
       Note that \eqref{aux} is trivial unless $t\ge w_m\log\left(2(m+1)\right)$ (as the right-hand side exceeds $1$). Therefore from now on we will consider only $t$ satisfying this lower bound. In particular, setting $p = 1/5$, we have $t \ge  \frac{2}{p}(2/\alpha)^\frac{1}{\alpha}c$ and $ t\ge 4^\frac{1}{\alpha} \frac{2c}{p} \log(n)^\frac{1}{\alpha}$.
       Using the union bound and the assumption $3)$, we get (denoting for brevity $\E_i\left(\cdot\right) = \cmv[\cdot][\ccF_i]$)
       \begin{equation*}
          \begin{split}
             \P\left(\sup_{1\le k \le n}\left|\sum_{i=1}^k X_i\right| > t \right) & \le \P\left(\sup_{1\le k \le n}\left|\sum_{i=1}^k Z_i\right| > t(1-p) \right)  + \P\left(\sup_{1\le i \le n}\left|\E_i X_{i+1} - \E_0 X_{1}\right| > tp \right)	\\
          \end{split}
       \end{equation*}
       \begin{equation}\label{B7}
          \qquad\le \P\left(\sup_{1\le k \le n}\left|\sum_{i=1}^k Z_i\right| > t(1-p) \right) + 2\P\left(\sup_{1\le i \le n}\left|\E_{i - 1}X_{i}\right| > \frac{tp}{2} \right).
       \end{equation}

       By another application of the union bound together with Lemma \ref{TICMVL} and stationarity of $\left(\E_{i - 1}X_i\right)_i$ we obtain
       $$2\P\left(\sup_{1\le i \le n}\left|\E_{i - 1}X_{i}\right| > \frac{tp}{2} \right) \le 2n \P\left(\left|\E_0 X_{1} \right| > \frac{tp}{2} \right) \le 12n\exp\left(-\frac{p^\alpha t^\alpha}{2(2c)^\alpha}\right).$$
       Notice that
       $$12n\exp\left(-\frac{p^\alpha t^\alpha}{2(2c)^\alpha}\right) = 12\left[n\exp\left(-\frac{p^\alpha t^\alpha}{4(2c)^\alpha}\right)\right]\exp\left(-\frac{p^\alpha t^\alpha}{4(2c)^\alpha}\right) \le 12\exp\left(-\frac{p^\alpha t^\alpha}{4(2c)^\alpha}\right),$$
       where the inequality is a consequence of  the estimate $ t\ge 4^\frac{1}{\alpha} \frac{2c}{p} \log(n)^\frac{1}{\alpha}$. It follows that
       \begin{equation} \label{IE_CE_SUP}
          2\P\left(\sup_{1\le i \le n}\left|\E_{i - 1} X_{i}\right| > \frac{pt}{2} \right) \le 12\exp\left(-\frac{p^\alpha t^\alpha}{4(2c)^\alpha}\right) = 12\exp\left(-\frac{t^\alpha}{4(10c)^\alpha}\right).
       \end{equation}
       In order to deal with $\P\left(\left|\sum_{i=1}^n Z_i\right| > t(1-p) \right)$ we start with splitting this sum into $m+1$ parts and using the union bound, namely
       \begin{equation*}
          \P\left(\sup_{1\le k \le n}\left|\sum_{i=1}^k Z_i\right| > t(1-p) \right) \le \sum_{j=0}^m \P\left(\sup_{1\le k \le n}\left|\sum_{1\le i \le k, m+1|i-j} Z_i\right| > \frac{t(1-p)}{m+1} \right).
       \end{equation*}
       Now, to each summand on the right-hand side of the above inequality we will apply the estimate for the independent case obtained at the beginning of this proof. Setting $M = (24\alpha^{-3} \log{n})^\frac{1}{ \alpha}c$ and taking into account \eqref{asymptotic_variance_formula} we obtain

       \begin{equation*}
          \begin{split}
             &\frac{1}{m + 1}\P\left(\sup_{1\le k \le n}\left|\sum_{i=1}^k Z_i\right| > t(1-p) \right) \le \frac{1}{m + 1} \sum_{j=0}^m \P\left(\sup_{1\le k \le n}\left|\sum_{1\le i \le k, m+1|i-j} Z_i\right| > \frac{t(1-p)}{m+1} \right)\\
             &\quad\le\exp(8)\exp\left(-\frac{t^\alpha}{\frac{16}{\alpha} (8(m+1)c)^\alpha}\right)+ 2\exp\left(-\frac{(1-p)^2 t^2 }{\frac{72}{25}(m+1)\left[\left(n+m+1\right)\sigma_\infty^2+\frac{8}{5}(1-p)tM\right]}\right)\\
          \end{split}
       \end{equation*}
       \begin{equation}\label{B10}
          \le \exp(8)\exp\left(-\frac{t^\alpha}{\frac{16}{\alpha} (8(m+1)c)^\alpha}\right)+ 2\exp\left(-\frac{ t^2 }{(m+1)\left[5\left(n+ m + 1\right)\sigma_\infty^2+2tM\right]}\right).
       \end{equation}
       Finally using \eqref{B7}, \eqref{IE_CE_SUP} and \eqref{B10} we get
       \begin{equation*}
          \begin{split}
             \P\left(\sup_{1\le k \le n}\left|\sum_{i=1}^k X_i\right| > t \right)
             \le& 12\exp\left(-\frac{t^\alpha}{4(10c)^\alpha}\right)+ (m+1)\exp(8)\exp\left(-\frac{t^\alpha}{\frac{16}{\alpha} (8(m+1)c)^\alpha}\right)	\\
             & \;  + 2(m+1)\exp\left(-\frac{t^2}{5(m+1)\left(n+m+1\right)\sigma_\infty^2+2(m+1)t M}\right).
          \end{split}
       \end{equation*}
       To conclude $\eqref{aux}$ it is now enough to note that the second summand on the right-hand side above dominates the first one.

       To finish the proof of the lemma it remains to show \eqref{asymptotic_variance_formula}. Firstly, we address the variance of $Z_i$, which can be easily calculated by using the properties of conditional expectation. We have (recall the notation $\E_i\left(\cdot\right) = \cmv[\cdot][\ccF_i]$)
       \begin{equation*}
          \begin{split}
             \E Z_i^2  = \E\big[ X_i^2+\E^2_i X_{i+1}+\E^2_{i - 1} X_{i} -2\E_{i}X_{i+1}\E_{i - 1} X_{i} - 2 X_i\E_{i - 1} X_{i}+2 X_i\E_{i} X_{i+1}\big].
          \end{split}
       \end{equation*}
       Since $\E X_i\E_{i - 1}X_{i} =\E \E^2_{i - 1}X_i$, $\E \E_{i} X_{i+1}\E_{i - 1} X_{i}=\E X_{i + 1} \E_{i - 1}X_{i}$ and $ X_i\E_{i} X_{i+1} = \E_{i}( X_i X_{i+1})$, we obtain
       \begin{equation*}
          \begin{split}
             \E Z_i^2 & = \E\Big( X_i^2+\E^2_{i}X_{i+1}-\E^2_{i - 1} X_{i} - 2 X_{i+1}\E_{i - 1} X_{i} + 2 X_i  X_{i+1} \Big)\\
             & = \E\left( X_i^2 + 2 X_i  X_{i+1}\right)-2\E \left(X_{i+1}\E_{i - 1} X_{i}\right)+  \E\left( \E^2_{i}X_{i+1}-\E^2_{i - 1}X_{i}\right). \\
          \end{split}
       \end{equation*}
       The variance formula in \eqref{asymptotic_variance_formula} follows by observing that due to $3)$, $\E\left( \E^2_{i}X_{i+1}-\E^2_{i - 1}X_{i}\right) = 0$, whereas by $4)$,
       $\E \left(X_{i+1}\E_{i - 1} X_{i}\right) = 0$.

       Now we will demonstrate the upper-bound on $ \|Z_i\|_{\psi_\alpha}$ in \eqref{asymptotic_variance_formula}. Using the triangle inequality (cf. Lemma \ref{triangle ineq alpha}) twice and then Lemma \ref{ONCMVL} we obtain
       \begin{equation} \label{Conditional_quasi_norm_bound}
          \begin{split}
             \|Z_i\|_{\psi_\alpha} & \le  2^{\frac{1}{\alpha} - 1}\|X_i\|_{\psi_\alpha} + 2^{\frac{1}{\alpha} - 1} \|\E_i X_{i+1} - \E_0 X_{1}\|_{\psi_\alpha} \le 2^{\frac{1}{\alpha}}\|X_i\|_{\psi_\alpha} + 2^{\frac{2}{\alpha}-1} \|\E_0 X_{1}\|_{\psi_\alpha} \\
             &\le 2^{\frac{1}{\alpha}}\|X_i\|_{\psi_\alpha} + 2^{\frac{2}{\alpha}-1} (2/\alpha)^\frac{1}{\alpha} \|X_{1}\|_{\psi_\alpha} \le \|X_{1}\|_{\psi_\alpha} \left(2^{\frac{1}{\alpha}}+\frac{1}{2} (8/\alpha)^\frac{1}{\alpha}\right) \le c (8/\alpha)^\frac{1}{\alpha}.
          \end{split}
       \end{equation}
       This concludes the proof of the lemma.
    \end{proof}

    \begin{rema}\em
       If $(X)_{i\ge 0}$ is a $1$-dependent, centered and stationary Markov chain such that $\|X_i\|_\infty \le M <\infty$ then the assumptions of the above lemma are satisfied with $m = 2$ and $\ccF_i = \sigma\left\{X_j \;|\; j\le i\right\}$.
       If  $(\xi_i)_{i\ge 0}$ are i.i.d. random variables and $f:\R^2\rightarrow \R$ is a bounded, Borel function such that $X_i = f(\xi_i, \xi_{i+1})$ are centered then we can take $\ccF_i = \sigma\{\xi_j \;|\; j\le i+1\}$ and notice that the assumptions of the above lemma are satisfied with $m = 1$.
    \end{rema}
    \begin{rema}\label{zero_variance}\em
       It is worth noticing that $\sigma_\infty^2$ may be equal to 0 in case of $1$-dependent processes $(X_i)_{i \in \N}$. Take for example $X_i= \xi_{i+1} - \xi_i$ where $(\xi_i)_{i \in \N}$ are i.i.d. random variables. It turns out (cf. \cite{SJ}) that the reverse is true, that is if for a $1$-dependent, bounded stationary process $(X_i)_{i \in \N}$ we have $\sigma_\infty^2 = 0$ then there exists an i.i.d. process $(\xi_i)_{i \in \N}$ such that $X_i= \xi_{i+1} - \xi_i$.
    \end{rema}
    \begin{lema}[Bernstein inequality for random sums]\label{1DBI_ST}
       Let $(X_i)_{i\geq 0}$ be a $1$-dependent sequence of centered random variables such that $\E\exp(c^{-\alpha}|X_i|^\alpha)\leq 2$ for some $\alpha \in (0,1]$ and $c \ge 1$. Moreover, let $N\le n\in\N$ be an $\N$-valued bounded random variable. Assume that we can find a filtration $\ccF = \left(\ccF_i\right)_{i\ge 0}$ such that for $Z_i =X_i + \E\left(X_{i+1}|\ccF_i\right) -\E\left(X_{i}|\ccF_{i-1}\right)$ we have the following:
       \begin{enumerate}
          \item[0)] $X_i$ is $\ccF_i$ measurable,
          \item[1)] $N$ is a stopping time with respect to $\ccF$,
          \item[2)] $(Z_i)_{i\ge 1}$ is stationary,
          \item[3)] For each $j \in\N$ process $(Z_i)_{i\ge j + 3}$ is independent of $\ccF_{j}$,
          \item[4)] $\left(\E\left(X_{i}|\ccF_{i-1}\right)\right)_{i\ge 1}$ is stationary,
          \item[5)] $ \E(X_i|\ccF_{i-1})$ is independent of $X_{i+1}$ for all $i \ge 1$.
       \end{enumerate}
       Then for any $t > 0$ and $a>0$,
       \begin{equation}\label{aux3}
          \P\left(\left|\sum_{i=1}^N X_i\right| > t \right)  \le 4\exp(8) \exp\left(-\frac{t^\alpha}{u c^\alpha}\right) + 9 \exp\left(-\frac{t^2}{v \sigma_\infty^2 + w t}\right),
       \end{equation}
       where $u = \frac{16 \cdot 26^\alpha}{\alpha}$, $v = 102 a$,  $w = 14M\max\left(2,  \sqrt{\|\left(\lceil N/3 \rceil - a + 1\right)_+ \|_{\psi_1}}\right)$ and \\
       $M = c(24\alpha^{-3} \log{n})^\frac{1}{ \alpha}$.
       \begin{proof}
          Observe that 0) and 3) imply $2$-dependence of the process $(Z_i)_{i\ge 1}$. Therefore the filtration $\mathcal{F}$ satisfies all the assumptions of Lemma \ref{1DUBI} and thus \eqref{asymptotic_variance_formula} holds. Note also that without loss of generality we may assume that $t\ge w\log 9 $ (otherwise the right-hand side of \eqref{aux3} is at least one). Fix $s = (8\sqrt{2}\log 9)^{-1}$. Using the union bound we get ($\E_i\left(\cdot\right) = \cmv[\cdot][\ccF_i]$)
          \begin{equation}\label{eq:first-union-bound}
             \P\left(\left|\sum_{i=1}^N X_i\right| > t \right) \le \P\left(\left|\sum_{i=1}^N Z_i\right| > t(1-s) \right) + 2\P\left(\sup_{1\le i \le n}\left| \E_{i - 1}X_i \right| > \frac{ts}{2} \right).
          \end{equation}
          Now using $\textnormal{Lemma \ref{TICMVL}, }$ $ts/2 \ge c\left(\frac{2}{\alpha}\right)^\frac{1}{\alpha}$, $t\ge w\log 9 $ and $n\exp\left(-\frac{(st)^\alpha}{4(2c)^\alpha}\right) \le 1$ we obtain
          \begin{equation}\label{final_0}
             \begin{split}
                2\P\left(\sup_{1\le i \le n}\left|\E_{i - 1}X_{i}\right| > \frac{st}{2} \right)  \le 2n\P\left(\left|\E_0 X_{1}\right| > \frac{st}{2} \right)   \le  12\exp\left(-\frac{(st)^\alpha}{4(2c)^\alpha}\right).
             \end{split}
          \end{equation}
          Next we take care of the other term on the right-hand side of \eqref{eq:first-union-bound}. Firstly we split the sum
          \begin{equation}\label{split_final}
             \P\left( \left|\sum_{i=1}^N Z_{i} \right| > t(1-s) \right)  \le \sum_{j=0}^2   \P\left( \left|\sum_{1\le i\le N,\; 3|(i+j)} Z_i \right| > \frac{t(1-s)}{3} \right).
          \end{equation}
          Now we will consider the $j$th summand of the above sum. Let us take $r = \frac{3}{8\sqrt{2}\log(9)}$ and notice that there is function $f_j:\N \rightarrow \N$ such that for any $n\in \N$, $\left\lfloor \frac{n}{3} \right\rfloor \le  f_j(n) \le \left\lceil \frac{n}{3} \right\rceil$ and
          \begin{equation}\label{aaaaaaaaaaaaaaaaaaaaa}
             \begin{split}
                & \P\left( \left|\sum_{1\le i\le N,\; 3|i + j} Z_{i} \right| > t(1-s)/3 \right) = \P\left( \left|\sum_{1\le i\le f_j(N)} Z_{3i - j} \right| > t(1-s)/3 \right)\\
                &\le \P\left( \left|\sum_{1\le i\le \lceil N/3 \rceil + 1} Z_{3i - j } \right| > t(1-r)(1-s)/3 \right) + \P\left(2\sup_{k\le n+6} \left|Z_{k} \right| >t r(1-s)/3\right).\\
             \end{split}
          \end{equation}
          Due to $\|Z_i\|_{\psi_\alpha} \le  c (8/\alpha)^\frac{1}{\alpha}$ (cf. \eqref{asymptotic_variance_formula}) and Lemma \ref{OCI} along with $t \ge w \log(9)$, $n \ge 2$  (for $n = 1$ the result of the lemma is trivial) we get
          \begin{equation} \label{final00}
             \begin{split}
                &\P\left(2\sup_{k\le n+6} \left|Z_{k} \right| >\frac{t r(1-s)}{3}\right)  \le (n + 6) 	\P\left(\left|Z_{k} \right| >\frac{t r(1-s)}{3}\right) \\
                & \le 2(n + 6) \exp\left(-\frac{\alpha(tr(1 - s))^\alpha}{8(3c)^\alpha}\right) \le 2\exp\left(-\frac{\alpha(tr(1 - s))^\alpha}{16(3c)^\alpha}\right).
             \end{split}
          \end{equation}
          To handle the first summand on the right-hand side of  \eqref{aaaaaaaaaaaaaaaaaaaaa} let us  fix $j$  and denote $\gamma_i :=  Z_{3i + 3 - j}$, $\ccG_i := \ccF_{3i - j}$, $T:= \lceil N/3 + 1\rceil  \le \lceil n/3 \rceil + 1$.
          Using the assumptions on the filtration $\mathcal{F}$ and \eqref{asymptotic_variance_formula} it is straightforward to check that the following properties hold:
          \begin{enumerate}
             \item  $\gamma_i$ are independent,
             \item $\E \gamma_i = 0$, $\E \gamma_i^2 = \sigma_\infty^2$, $\|\gamma_i\|_{\psi_\alpha} \le c (8/\alpha)^\frac{1}{\alpha}$,
             \item  $\gamma_{i-1}$ is $\ccG_i$ measurable,
             \item  $\gamma_i$ is independent of $\ccG_i$,
             \item  $T$ is a stopping time with respect to the filtration $\ccG_i$.
          \end{enumerate}
          This is precisely the setting of Proposition $4.4$. ii) from \cite{AB} which applied with $\epsilon := 1$, $p:= \frac{\sqrt{2}}{\sqrt{2}-1}$ and $q := \sqrt{2}$ gives
          that for any $a > 0$,
          \begin{equation}\label{final}
             \begin{split}
                &\P\left( \left|\sum_{1\le i\le \lceil N/3 \rceil + 1} Z_{3i - j } \right| > t(1-r)(1-s)/3 \right)\\
                & \le \exp(8)\exp\left(-\frac{(t(1-r)(1-s))^\alpha}{2(3(2+\sqrt{2})\hat{c})^\alpha}\right)+ 3\exp\left(-\frac{(t(1-r)(1-s))^2}{72 a \sigma_\infty^2 +6\sqrt{2}\mu (1-r)(1-s)t}\right),
             \end{split}
          \end{equation}
          where
          $$\mu = \max\left(\frac{8M}{3}, 2\sigma_\infty \sqrt{\|\left(\lceil N/3 \rceil - a + 1\right)_+ \|_{\psi_1}}\right), \qquad \hat{c} =  c \left(\frac{8}{\alpha}\right)^\frac{1}{\alpha}.$$

          Using \eqref{asymptotic_variance_formula}, Lemma \ref{moment_estimation} with $Y = \frac{\alpha Z^\alpha}{8c^\alpha}$ and $\beta = \frac{2}{\alpha}$, together with the gamma function estimate  $\Gamma(x) \le \left(\frac{x}{2}\right)^{x - 1} \textnormal{ for } x \ge 2$ (see Theorem $1$ in \cite{LC}) we get
          \begin{equation*}
             \sigma_\infty^2 = \E Z_1^2   \le 2c^2\left(\frac{8}{\alpha}\right)^\frac{2}{\alpha}\Gamma\left(\frac{2}{\alpha} + 1\right) \le \frac{4}{\alpha}c^2\left(\frac{8}{\alpha}\right)^\frac{2}{\alpha}\Gamma\left(\frac{2}{\alpha}\right)\le  4 c^2\left(\frac{8}{\alpha^2}\right)^{\frac{2}{\alpha}},
          \end{equation*}
          which implies that $\sigma_\infty \le \frac{2}{3}M$ and as a consequence,
          \begin{equation*}
             \mu \le \frac{4}{3}M b, \qquad  b = \max\left(2,  \sqrt{\|\left(\lceil N/3 \rceil - a + 1\right)_+ \|_{\psi_1}}\right).
          \end{equation*}
          Therefore \eqref{final} reduces to
          \begin{equation*}
             \begin{split}
                &\P\left( \left|\sum_{1\le i\le \lceil N/3 \rceil + 1} Z_{3i - j } \right| > t(1-r)(1-s)/3 \right)\\
                & \le \exp(8)\exp\left(-\frac{(t(1-r)(1-s))^\alpha}{2(3(2+\sqrt{2})\hat{c})^\alpha}\right)+ 3\exp\left(-\frac{(t(1-r)(1-s))^2}{72 a \sigma_\infty^2 +8\sqrt{2}M b(1-r)(1-s)t}\right).\\
             \end{split}
          \end{equation*}
          Combining the above inequality with \eqref{eq:first-union-bound}--\eqref{final00} we obtain
          \begin{equation*}
             \begin{split}
                &\P\left(\left|\sum_{i=1}^N X_i\right| > t \right)  \le   12\exp\left(-\frac{(st)^\alpha}{4(2c)^\alpha}\right) + 6\exp\left(-\frac{\alpha(tr(1 - s))^\alpha}{16(3c)^\alpha}\right)\\
                &\quad + 3\exp(8)\exp\left(-\frac{(t(1-r)(1-s))^\alpha}{2(3(2+\sqrt{2})\hat{c})^\alpha}\right) +9\exp\left(-\frac{(t(1-r)(1-s))^2}{72 a \sigma_\infty^2 +8\sqrt{2}M b(1-r)(1-s)t}\right).
             \end{split}
          \end{equation*}
          To conclude it is now enough to recall that $r = 3(8\sqrt{2}\log(9))^{-1}$, $s = (8\sqrt{2}\log 9)^{-1}$ and do some elementary calculations.
       \end{proof}
    \end{lema}
\section{Bounds on the number of regenerations}
    We will now obtain a bound on the stopping time $N$, introduced in \eqref{eq:definition-of-Theta-N}.
    To this end we will use the $\psi_1$ version of Bernstein inequality, which follows easily from the classical moment version of this inequality (see, e.g., Lemma 2.2.11 in \cite{VaartWellnerPsi1}), by observing  that for $k \ge 2$, $\E |\xi|^k \le k! \|\xi\|_{\psi_1}^k = k!M^{k-2}v/2$, where $M = \|\xi\|_{\psi_1}$, $v = 2\|\xi\|_{\psi_1}^2$.
    \begin{lema} ($\psi_1$ Bernstein's inequality.) \label{psi1 BE}
       If $(\xi_i)_i$ is a sequence of independent centered random variables such that $\sup_i \|\xi_i \|_{\psi_1}\leq \tau$, then
       $$\P\left(\sum_{i=1}^n \xi_i \geq t\right)\leq \exp\left(- \frac{t^2}{4n\tau^2+2\tau t}\right).$$
    \end{lema}
    \begin{lema}\label{Regeneration}
       If $\|\tau_1-\tau_0\|_{\psi_1} \leq d $ then for any $p>0$,
       \begin{equation}\label{eq:bound-on-N}
          \P\left(N > \left\lceil{(1+p)n\left[\E(\tau_1-\tau_0)\right]^{-1}}\right\rceil \right)\leq \exp(1)\exp\left(-\frac{p n \E(\tau_1-\tau_0)}{K_p d^2 }\right),
       \end{equation}
       where $K_p = L_p + 16/L_p$ and $L_p = \frac{16}{p} + 20$. Moreover,  the function $p\rightarrow K_p$ is decreasing on $\R_+$ (in particular $K_p \ge K_\infty = \frac{104}{5}$) and if $p=2/3$ then $\frac{1}{p} K_{p} \le 67$ .
       \begin{proof}
          For convenience, let $T_i = \tau_i - \tau_{i - 1}$ for $i \ge 1$.
          Firstly, notice that without loss of generality we may assume that $np \ge L_p \E T_1$. Indeed, otherwise, using $\E T_1 \le d$ we obtain
          \begin{equation*}
             \exp(1)\exp\left(-\frac{p n \E T_1}{K_p d^2 }\right) \ge \exp(1)\exp\left(-\frac{L_p \E^2T_1}{K_p d^2 }\right) \ge \exp\left(1 -\frac{L_p}{K_p}\right)  \ge 1.
          \end{equation*}
          Thus, from now on we consider $n$ such that $np \ge L_p \E T_1$. For $A = (1+p)n\left[\E T_1\right]^{-1} \ge 1$ we get
          \begin{equation}
             \begin{split}
                & \P(N > \lceil{A}\rceil ) \le \P(\tau_{\lceil{A}\rceil}-\tau_0\leq n) \le \P\left( \sum_{i=0}^{\lceil{A}\rceil -1}T_{i + 1} -\E T_{i + 1} \le  n -A \E T_{1}\right) \\
                & = \P\left(\sum_{i=0}^{\lceil{A}\rceil - 1 }T_{i + 1} -\E T_{i + 1} \le n - (1+p)n\right) = \P\left(\sum_{i=0}^{\lceil{A}\rceil -1}T_{i + 1} -\E T_{i + 1} \le -np\right).
             \end{split}
          \end{equation}
          Now we have $\|T_{i+1} -\E T_{i+1}\|_{\psi_1} \leq 2d$ so using Lemma \ref{psi1 BE},  $\E T_1 \leq d$  and $np \ge L_p \E T_1$ we get
          \begin{equation*}
             \begin{split}
                &\P\left(N > \left\lceil(1+p)n \left[\E T_1\right]^{-1} \right\rceil\right) \leq  \exp\left(-\frac{p^2 n^2}{4(A+1) 4d^2 + 4d np}\right) \\
                & =  \exp\left(-\frac{p^2 n^2}{16 d^2 \left[ (1+p)n\left[\E T_1\right]^{-1} + 1\right] + 4d np}\right) =  \exp\left(-\frac{p n \E T_{1}}{16 d^2 \left( \frac{1+p}{p} +\frac{\E T_{1}}{pn}\right) + 4d \E T_{1}}\right)\\
                & \le    \exp\left(-\frac{p n \E T_{1}}{16 d^2 \left( \frac{1+p}{p} +\frac{1}{L_p}\right) + 4d^2}\right) = \exp\left(-\frac{p n \E T_{1}}{K_p d^2 }\right) \le \exp\left(1 - \frac{p n \E T_{1}}{K_p d^2 }\right),
             \end{split}
          \end{equation*}
          which finishes the proof of \eqref{eq:bound-on-N}. The properties of $K_p$ follow from easy computations.
       \end{proof}
    \end{lema}

    The following lemma is a standard consequence of the tail estimates given in Lemma \ref{Regeneration}. Its proof, based on integration by parts, is analogous to that of  Lemma $5.4$ in \cite{AB} and is therefore omitted.

    \begin{lema} \label{psi_1_exp_N}
       Suppose that $\|\tau_1-\tau_0\|_{\psi_1} \leq d $ for some $d>0$. Then for any $p>0$,
       \begin{equation*}
          \left\|\left(N - a\right)_+\right\|_{\psi_1} \le \frac{4 K_p d^2}{\left[\E(\tau_1-\tau_0)\right]^2}\le \frac{4 K_p d^2}{m^2},
       \end{equation*}
       where  $a = (1+p) n \left[\E(\tau_1-\tau_0)\right]^{-1}$, $K_p = L_p + \frac{16}{L_p}$ and $L_p =  \frac{16}{p} + 20$.  Moreover,
       \begin{equation*}
          \frac{d^2 K_p}{\left[\E(\tau_1-\tau_0)\right]^2} \ge K_p \ge K_\infty.
       \end{equation*}
    \end{lema}

\section{Proofs of Theorems \ref{BI}, \ref{BI2} and \ref{SBI}}
    In this section we will prove our main results. The structure of proofs of Theorems \ref{BI} and \ref{BI2} is similar, and they contain a common part, which we will present first in Sections \ref{secH} and \ref{sec:T}. The proof of Theorem \ref{BI} will be concluded in Section \ref{sec:BI-proof} and the proof of Theorem \ref{BI2} in Section \ref{sec:BI2-proof}. Theorem \ref{SBI} will be obtained as a corollary to Theorem \ref{BI} in Section \ref{sec:SBI}.

    Let us thus pass to the proofs of Theorems \ref{BI} and \ref{BI2}. Assume that $m|n$. The argument will be based on the approach of \cite{A} and \cite{AB} (see also \cite{C} and \cite{DGM}) and will rely on the decomposition
    \begin{equation} \label{abs_block_decomposition}
       \left|\sum_{i=0}^{n-1}f(\Upsilon_i)\right|\leq H_n + M_n  + T_n,
    \end{equation}
    where
    \begin{equation}
       \begin{split}
          H_n  & = \left|\sum_{i=0}^{\tau_0 /m} \Theta_i \1_{N>0}  + \1_{N=0} \sum_{i=0}^{n/m - 1} \Theta_i \right|,\; M_n = \left|\sum_{i=1}^{N} \chi_{i-1}(f)\right|,\\
          T_n &=\left|\1_{N>0} \sum_{k=n}^{\tau_N + m - 1} f(\Upsilon_k)\right|, \; N = \inf\{i\geq 0 \;|\; \tau_i +  m - 1\geq n-1\}.
       \end{split}
    \end{equation}
    The proof will be divided into three main steps. In the first two (common for both theorems) we will get easy bounds on tails of $H_n$ and $T_n$. The main, third step will be devoted to obtaining two different estimates on the tail of $M_n$. To this end we will use Lemmas \ref{1DUBI}, \ref{Regeneration} (for the proof of Theorem \ref{BI}) and Lemmas \ref{1DBI_ST}, \ref{psi_1_exp_N} (for Theorem \ref{BI2}).

    \subsection{Estimate on \texorpdfstring{$H_n$}{}}\label{secH}
        Using $\{N=0\}\subset \{\tau_0 \ge n - m\}$, the definition of $\a$ (see \eqref{eq1}) and Lemma \ref{OCI} we get
        \begin{equation}\label{ine_first_block}
           \begin{split}
              \P_{x^*}(H_n>t) & \le \P_{x^*}\left(\1_{N>0}\sum_{i=0}^{\tau_0 /m}\left| \Theta_i \right|   + \1_{N=0} \sum_{i=0}^{n/m - 1} \left| \Theta_i \right|  > t\right)\le \P_{x^*}\left(\sum_{i=0}^{\tau_0 /m}\left| \Theta_i \right|  > t\right)\\
              &\le 2 \exp\left(-\frac{t^\alpha}{\a^\alpha}\right).
           \end{split}
        \end{equation}

    \subsection{Estimate on \texorpdfstring{$T_n$}{}}\label{sec:T}
        By repeating verbatim the easy argument presented in the proof of Theorem 5.1 in \cite{AB}, we obtain
        \begin{equation}\label{ine_tail_block}
           \P\left( \left|T_n\right| > t \right)\leq 2\left[\delta\pi(C)\right]^{-1} \exp\left(-\frac{t^\alpha}{\b^\alpha}\right).
        \end{equation}
        We skip the details.
    \subsection{Proof of Theorem \ref{BI}}\label{sec:BI-proof}
        Recall that $M = \c(24\alpha^{-3} \log{n})^\frac{1}{ \alpha}$ and note that without loss of generality we can assume that $t \ge 8 M \log 6$. Otherwise \eqref{eq:first-theorem} is trivial as the right hand side is greater than or equal to $1$.
        Fix $p = 2/3$. We have ($A := \left\lceil{(p+1)n(\E(\tau_1-\tau_0))^{-1}}\right\rceil$)
        \begin{equation}\label{eq:nie-mam-pomyslu-na-nazwe}
           \begin{split}
              &\P\left(M_n \ge t \right) = \P\left(M_n \ge t ,\; N\le A \right) + \P\left(M_n \ge t, N > A\right)\\
              & \le \P\left( \sup_{1\le k \le  A} \left|\sum_{i=1}^{k} \chi_{i-1}\right| \ge t  \right) + \P\left(N > A\right).
           \end{split}
        \end{equation}

        To control the first summand on the right-hand side of the above inequality we will apply Lemma \ref{1DUBI} with $m=2$, $X_i := \chi_{i} = F(\Xi_{i + 1})$ (cf. \eqref{siconditional}),  $c := \c$ and $n := A$. Assuming that the assumptions of the lemma are satisfied (we will verify them later on), we obtain (in the first line we use stationarity of $(\Xi_i)_{i\ge 1}$),
        \begin{equation*}
           \begin{split}
              & P := \P\left( \sup_{1\le k \le A} \left|\sum_{i=1}^{k} \chi_{i-1}\right| \ge t  \right) = \P\left( \sup_{1\le k \le  A} \left|\sum_{i=1}^{k} F(\Xi_{i+1})\right| \ge t  \right) \\
              & \le  6\exp(8)\exp\left(-\frac{t^\alpha}{\frac{16}{\alpha} (24 \c)^\alpha}\right)
              + 6\exp\left(-\frac{t^2}{15 \left(\left \lceil{(p+1)n(\E(\tau_1-\tau_0))^{-1}}\right \rceil+3\right)\sigma_\infty^2+6 t M}\right)\\
           \end{split}
        \end{equation*}
        \begin{equation}\label{eq:auxiauxi}
           \;\;\;\le  6\exp(8)\exp\left(-\frac{t^\alpha}{\frac{16}{\alpha} (24 \c)^\alpha}\right)
           + 6\exp\left(-\frac{t^2}{15 \left((p+1)n(\E(\tau_1-\tau_0))^{-1}+4\right)\sigma_\infty^2+6 t M}\right)
        \end{equation}

        Recall that by \eqref{asymp_variances_formula}, $\sigma_\infty^2 = \sigma_{Mrv}^2 \E (\tau_1-\tau_0)$. We will now obtain a comparison between $\sigma^2_\infty$ and $tM$, which will allow us to reduce the above estimate to one in which the subgaussian coefficient is expressed only in terms of $\sigma_{Mrv}^2$.  Thanks to Lemma \ref{moment_estimation} applied with $\Lambda := (\chi_1/\c)^\alpha$ and $\beta := 2/\alpha$ we have
        $$\sigma_\infty^2 \le 3\E\chi^2_1 \le 3\c^2 \Gamma(2/\alpha + 1) \le 3\c^2 (2/\alpha)^{\frac{2}{\alpha} + 1},$$
        where the last inequality is a consequence of equation $4$ in \cite{LC}.
        Moreover, recalling the definition of $M$ and using the assumption $t \ge 8\log(6)M$, we obtain
        \begin{equation*}
           \begin{split}
              tM & \ge 8 \log(6) M^2 =8\log(6) \c^2(24 \alpha^{-3} \log(n))^\frac{2}{\alpha}\ge 16\cdot 8\log(6) \c^2(2/\alpha)^{\frac{2}{\alpha} + 1} \ge 76 \sigma_\infty^2.
           \end{split}
        \end{equation*}
        The last inequality in combination with \eqref{eq:auxiauxi} yields
        \begin{equation}\label{n}
           \begin{split}
              P &\le  6\exp(8)\exp\left(-\frac{t^\alpha}{\frac{16}{\alpha} (24 \c)^\alpha}\right)
              + 6\exp\left(-\frac{t^2}{15 (p+1)n\sigma_{Mrv}^2+7 t M}\right)\\
              & \le  6\exp(8)\exp\left(-\frac{t^\alpha}{\frac{16}{\alpha} (24 \c)^\alpha}\right)
              + 6\exp\left(-\frac{t^2}{25 n\sigma_{Mrv}^2+ 7 t M}\right).
           \end{split}
        \end{equation}

        In order to justify the above inequality it remains to verify the assumptions of Lemma \ref{1DUBI}. To this end take
        \begin{equation*}
           \ccF_i = \sigma\{\Xi_j\;|\; j\leq i+1 \}, \quad Z_i =\chi_i+\E(\chi_{i+1}|\ccF_{i})-\E(\chi_{i}|\ccF_{i-1}).
        \end{equation*}
        We will now strongly rely on the properties of the stationary sequence of one-dependent blocks $(\Xi_i)_{i\geq 1}$ stated in Remark \ref{rem:Rndom-blocks} together with \eqref{Excusions Definition} and \eqref{siconditional}. Since $\chi_i = F(\Xi_{i+1})$, the assumption 0) of Lemma \ref{1DUBI} is trivially satisfied. To prove 2), observe that $\E(\chi_{i+1}|\ccF_{i}) = \E(\chi_{i+1}|\Xi_{i+1}) = G(\Xi_{i+1})$ for some measurable function $G$, and so the sequence $\left(Z_i\right)_{i\ge 1}=\left(F(\Xi_{i+1})+G(\Xi_{i+1}) - G(\Xi_{i})\right)_{i\ge 1}$ is stationary (as a function of the stationary sequence $\left(\Xi_i\right)_{i\ge 1}$). The sequence $\left(\Xi_i\right))_{i\ge 1}$ is $1$-dependent, which clearly implies that $\left(Z_i\right)_{i\ge 1}$ is $2$-dependent, i.e., the assumption 2) of the lemma. The assumption 3), i.e., the stationarity of the sequence $\left(\E(\chi_i|\ccF_{i-1})\right)_{i\ge 1} = \left(G(\Xi_i)\right)_{i\ge 1}$ follows again by stationarity of $\left(\Xi_i\right)_{i\ge 1}$. Finally, using once more the fact that $(\Xi_i)_{i\ge 0}$ is one-dependent, we obtain that for any $i\ge 1$ the random variable $\E(\chi_i|\ccF_{i-1}) = G(\Xi_i)$ is independent of $\chi_{i+1} = F(\Xi_{i+2})$, which ends the verification of the assumptions of Lemma \ref{1DUBI} and proves \eqref{n}.

        Thus, in order to get a bound on $\P(M_n > t)$ it suffices to estimate the second term on the right-hand side of \eqref{eq:nie-mam-pomyslu-na-nazwe}. To this aim we use Lemma \ref{Regeneration} with $p =2/3$ and $d= \d$ obtaining
        \begin{equation*}
           \P\left(N > \left\lceil{(1+p)n(\E(\tau_1-\tau_0))^{-1}}\right\rceil \right)\leq \exp(1)\exp\left(-\frac{ n \E(\tau_1-\tau_0)}{67 \d^2 }\right).
        \end{equation*}
        In combination with \eqref{eq:nie-mam-pomyslu-na-nazwe} and \eqref{n} this gives
        \begin{equation*}
           \begin{split}
              \P\left(M_n \ge t \right) & \le	6\exp(8)\exp\left(-\frac{t^\alpha}{\frac{16}{\alpha} (24 \c)^\alpha}\right)
              + 6\exp\left(-\frac{t^2}{25 n\sigma_{Mrv}^2+ 7 t M}\right) \\
              &\; +\exp(1)\exp\left(-\frac{ n \E(\tau_1-\tau_0)}{67 \d^2 }\right).
           \end{split}
        \end{equation*}
        \noindent Combining the above inequality with \eqref{ine_first_block} and \eqref{ine_tail_block} we get
        \begin{equation*}
           \begin{split}
              &\P_x \left(\left|\sum_{i=0}^{n-1} f(\Upsilon_i)\right|>t\right)  \le \P\left(H_n\geq \frac{1-\sqrt{5/6}}{2}t\right)+ \P\left(M_n\geq \sqrt{5/6}t\right) + \P\left(T_n\geq \frac{1-\sqrt{5/6}}{2}t\right) \\
              & \le 2\exp\left(-\frac{t^\alpha}{(23\a)^\alpha}\right) + 2\left[\delta\pi(C)\right]^{-1} \exp\left(-\frac{t^\alpha}{(23\b)^\alpha}\right)  +\exp(1)\exp\left(-\frac{ n \E(\tau_1-\tau)}{67 \d^2 }\right)\\
              &\;\; +	6\exp(8)\exp\left(-\frac{t^\alpha}{\frac{16}{\alpha} (27 \c)^\alpha}\right)
              + 6\exp\left(-\frac{t^2}{30 n\sigma_{Mrv}^2+8 t M}\right).
           \end{split}
        \end{equation*}
        In order to finish the proof of Theorem \ref{BI} it is enough to substitute $\E(\tau_1-\tau_0) = \delta^{-1} \pi(C)^{-1}m$.

    \subsection{Proof of Theorem \ref{BI2}}\label{sec:BI2-proof}
        Recall that $M = \c(24\alpha^{-3} \log{n})^\frac{1}{ \alpha}$ and let $p>0$ be a parameter which will be fixed later on. We are going to apply Lemma \ref{1DBI_ST} with $X_i := \chi_{i} = F(\Xi_{i + 1})$, $c := \c$, $\ccF_i := \sigma\{\Xi_j\;|\; 0\le j \le i + 1\}$. Clearly $N$ is a stopping time with respect to $\ccF$. The remaining assumptions of Lemm \ref{1DUBI} can be verified in the same manner as in the proof of Theorem \ref{BI}. Let $a = (1+p)\frac{n}{3} \left[\E(\tau_1-\tau_0)\right]^{-1}$. By Lemma \ref{psi_1_exp_N} we get
        \begin{equation*}
           \begin{split}
              \left\| \left(\lceil N/3 \rceil - a + 1 \right)_+\right\|_{\psi_1} &\le \frac{1}{3}\left\| \left( N  - (1+p)n(\E(\tau_1-\tau_0))^{-1}  \right)_+\right\|_{\psi_1} + \frac{2}{\log 2}\\
              & \le \frac{4}{3}  \d^2 K_p +  \frac{2}{\log 2} \le \left(\frac{4}{3} +\frac{7}{50}\right) \d^2 K_p,
           \end{split}
        \end{equation*}
        where the last inequality follows from (recall the definition of $K_\infty$ from Lemma \ref{Regeneration})
        $$\frac{7}{50}K_p \ge \frac{7}{50}K_\infty = \frac{7}{50}\cdot\frac{104}{5} \ge \frac{2}{\log 2}.$$
        Therefore $\max\left(2,  \sqrt{\|\left(\lceil N/3 \rceil - a + 1\right)_+ \|_{\psi_1}}\right) \le \sqrt{4/3 + 7/50} \sqrt{K_p} \cdot \d$ and we get that for arbitrary $p > 0$,
        \begin{equation*}
           \P\left( \left|\sum_{i=1}^N \chi_{i-1}\right|>t\right) \le 4\exp(8)\exp\left(-\frac{t^\alpha}{\frac{16}{\alpha}(26\c)^\alpha}\right) +9\exp\left(-\frac{t^2}{34(1+p) \sigma_{Mrv}^2 + 17 M \d t\sqrt{K_p}}\right).
        \end{equation*}
        \noindent Using the above inequality together with \eqref{ine_first_block}, \eqref{ine_tail_block} we obtain
        \begin{equation*}
           \begin{split}
              &\P_x \left(\left|\sum_{i=0}^{n-1} f(\Upsilon_i)\right|>t\right) \le \P\left(H_n\geq \frac{t}{54}\right)+ \P\left(M_n\geq \frac{26t}{27}\right) + \P\left(T_n\geq \frac{t}{54}\right) \\
              & \le 2\exp\left(-\frac{t^\alpha}{(54\a)^\alpha}\right) + 2\left[\delta\pi(C)\right]^{-1} \exp\left(-\frac{t^\alpha}{(54\b)^\alpha}\right) + 4\exp(8)\exp\left(-\frac{t^\alpha}{\frac{16}{\alpha}(27\c)^\alpha}\right)\\
              &\;\;  +9\exp\left(-\frac{t^2}{37(1+p) \sigma_{Mrv}^2 + 18 M \d t\sqrt{K_p}}\right),
           \end{split}
        \end{equation*}
        which concludes the proof of Theorem \ref{BI2}.

    \subsection{Proof of Theorem \ref{SBI}.}\label{sec:SBI}
        Denote $M = \norm[f][\infty]$ and notice that for $t > n M$ the left-hand side of \eqref{eq:third theorem} vanishes, so we may assume that $t \le n M$.
        Using \eqref{c bound} one can easily see that if $m|n$ then Theorem \ref{BI} applied with $\alpha = 1$ implies that
        \begin{equation}\label{simplify}
           \begin{split}
              &\P_x \left(\left|\sum_{i=0}^{n-1} f(\Upsilon_i)\right|>t\right) \le \left(2 +  2\left[\delta\pi(C)\right]^{-1}\right)\exp\left(-\frac{t}{46DM}\right)  + 6\exp(8)\exp\left(-\frac{t}{432DM}\right)\\
              &\;\;\qquad \qquad  + 6\exp\left(-\frac{t^2}{30 n\sigma_{Mrv}^2+ 192 t DM}\right) +\exp(1)\exp\left(-\frac{ n m}{67 \delta \pi(C)D^2 }\right).
           \end{split}
        \end{equation}
        The assumption $t \le nM$ yields
        $$\exp\left(-\frac{ n m}{67 \delta \pi(C)D^2 }\right) \le \exp\left(-\frac{ t m}{67 \delta \pi(C)MD^2 }\right),$$
        which plugged into \eqref{simplify} gives after some elementary calculations, that (recall $K =\exp(10) +  2\left[\delta\pi(C)\right]^{-1}$)
        \begin{equation}\label{m divides n}
           \begin{split}
              \P_x \left(\left|\sum_{i=0}^{n-1} f(\Upsilon_i)\right|>t\right) & \le K\exp\left(-\frac{t^2}{30n\sigma_{Mrv}^2 + 432t D^2 M \delta \pi(C) \log n}\right),
           \end{split}
        \end{equation}
        proving the theorem in the special case $m|n$.

        Now we consider the case $m\nmid n$. Define $\lceil n \rceil_m$ to be the smallest integer greater or equal to n, which is divisible by $m$. Notice that without loss of generality we can assume that $t > 4330D^2M\delta\pi(C)$ (otherwise the assertion of the theorem is trivial as the right-hand side of \eqref{eq:third theorem} exceeds one). Since $D^2\delta\pi(C) > m$ (recall $\E(\tau_1-\tau_0) = \delta^{-1} \pi(C)^{-1}m$), this implies that $t \ge 4330 Mm$. Moreover, as $t \le nM$, we also obtain that $n\ge 4330m$.

        Thus, for $p = 1/4330$ we have $\left|\sum_{i=n}^{\lceil n \rceil_m} f(\Upsilon_i)\right|\le Mm \le pt$, and as a consequence
        \begin{equation}
           \begin{split}
              \P_x \left(\left|\sum_{i=0}^{n-1} f(\Upsilon_i)\right|>t\right)  \le \P_x \left(\left|\sum_{i=0}^{\lceil n \rceil_m -1} f(\Upsilon_i)\right|>(1-p)t\right).
           \end{split}
        \end{equation}
        Now using \eqref{m divides n} and the inequality $n > 4330m$ we get 
        \begin{equation*}
           \begin{split}
              & \P_x \left(\left|\sum_{i=0}^{n-1} f(\Upsilon_i)\right|>t\right)  \le K\exp\left(-\frac{t^2}{31\lceil n \rceil_m \sigma_{Mrv}^2 + 433t D^2 M\delta \pi(C)\log n}\right)\\
              & \qquad \le K\exp\left(-\frac{t^2}{31 (n + m) \sigma_{Mrv}^2 + 433 t D^2 M\delta \pi(C)\log n}\right)\\
              & \qquad \le K\exp\left(-\frac{t^2}{32n \sigma_{Mrv}^2 +   433t D^2 M\delta\pi(C)\log n}\right). \\
           \end{split}
        \end{equation*}
        This concludes the proof of Theorem \ref{SBI}.

        \appendix

\section{Orlicz exponential norm}\label{appendix exponential norm}

    At the beginning recall the definition of the exponential Orlicz quasi-norm \eqref{definition of the exponential Orlicz norm} and note that if $\alpha \ge 1$ then $\| \cdot \|_{\psi_\alpha}$ is a norm whereas for $0< \alpha < 1$, $\| \cdot \|_{\psi_\alpha}$ is only a quasi-norm. More precisely, we have the following version of the triangle inequality (see Lemma $3.7$ in \cite{AB}).

    \begin{lema}[Triangle inequality for $\alpha \le 1$] \label{triangle ineq alpha}
       Fix $0 < \alpha \le 1$. Then for any random variables $X$, $Y$ we have
       $$\| X + Y \|_{\psi_\alpha} \le \left(\| X \|_{\psi_\alpha}^\alpha + \| Y \|_{\psi_\alpha}^\alpha\right)^{1/\alpha}  \le 2^{1/\alpha - 1}\left(\| X \|_{\psi_\alpha} + \| Y \|_{\psi_\alpha}\right).$$
    \end{lema}

    \noindent Now we present a moment estimation for random variables with bounded exponential moment.
    \begin{lema}\label{moment_estimation}
       If $Y$ is non negative random variable such that $\E \exp(Y) \le 2$ then for any $\beta > 0$ we have
       $$\E Y^\beta \le 2\Gamma(\beta + 1).$$
       Furthermore, if $\beta \in \N$ then one can replace the constant $2$ with $1$.
       \begin{proof} \normalfont
          If $\beta$ is a natural number then the claim follows from Taylor's expansion of $\exp(x)$. The general case is obtained by Markov's inequality, namely
          \begin{equation}
             \begin{split}
                \E Y^\beta & = \int_{0}^{\infty} \P\left(Y^\beta \ge t\right) dt = \int_{0}^{\infty} \P\left(\exp(Y) \ge \exp\left(t^\frac{1}{\beta}\right)\right) dt \\\
                & \le \int_{0}^{\infty} \frac{2}{\exp\left(t^\frac{1}{\beta}\right)} dt=\int_{0}^{\infty} \exp(-s)2\beta s^{\beta-1}ds = 2\beta\Gamma(\beta).
             \end{split}
          \end{equation}
       \end{proof}
    \end{lema}

    \noindent The next lemma allows to pass from the $\psi_\alpha$-norm of a random variable to the norm of its conditional expectation.

    \begin{lema}[Orlicz's norm of Conditional Mean Value]\label{ONCMVL}
       Let $0 < \alpha \le 1$. Assume that a random variable $X$ satisfies $\| X \|_{\psi_\alpha} < \infty$. Moreover, let $\ccF$ be some sigma field. Then
       $$\| \E(X|\ccF) \|_{\psi_\alpha} \le \left(1 + \frac{\log\left(\alpha \exp\left(\frac{1-\alpha}{\alpha}\right)\right)}{\log(2)}\right)^\frac{1}{\alpha}\| X \|_{\psi_\alpha} \le \left(\frac{2}{\alpha}\right)^\frac{1}{\alpha}\| X \|_{\psi_\alpha}.$$

       \begin{proof}
          Set $\phi_\alpha (x)= \exp(x^\alpha)$ for $x \ge 0$ and notice that $\phi_\alpha$ is concave on $(0, x_\alpha)$ and convex on $(x_\alpha, \infty)$, where  $x_\alpha = \left(\frac{1-\alpha}{\alpha}\right)^{1/\alpha}$. Define  $\Psi_\alpha$ to be a smallest convex function bigger or equal to $\phi_\alpha$ which is equal to $\phi_\alpha$ on $(x_\alpha, \infty)$, that is
          \[
             \Psi_\alpha(x)=
             \begin{dcases}
                \alpha \exp\left(\frac{1-\alpha}{\alpha}\right)(x x_\alpha^{\alpha - 1} + 1), & \text{if } 0\le x\le x_\alpha, \\
                \phi_\alpha(x),                                                               & \text{if }x\ge x_\alpha.
             \end{dcases}
          \]
          Clearly $\Psi_\alpha$ is a convex function on $\R_+$ and it is easy to see that $\phi_\alpha \le \Psi_\alpha \le \alpha \exp\left(\frac{1-\alpha}{\alpha}\right) \phi_\alpha$.
          Using these properties, Jensen's inequality and the definition of the Orlicz norm we get
          \begin{equation*}
                \E\phi_\alpha\left(\frac{\left|\E\left(X|\ccF \right)\right|}{\| X \|_{\psi_\alpha}}\right) \le  \E\Psi_\alpha\left(\frac{\left|\E\left(X|\ccF \right)\right|}{\| X \|_{\psi_\alpha}}\right)\le \E\Psi_\alpha\left(\frac{\left|X \right|}{\| X \|_{\psi_\alpha}}\right) \le 2\alpha \exp\left(\frac{1-\alpha}{\alpha}\right).
          \end{equation*}
          Put $c_\alpha = \left(1 + \frac{\log\left(\alpha \exp\left(\frac{1-\alpha}{\alpha}\right)\right)}{\log(2)}\right)^\frac{1}{\alpha} \ge 1$ and note that due to Jensen's inequality
          \begin{equation*}
	          \E\phi_\alpha\left(\frac{\left|\E\left(X|\ccF \right)\right|}{c_\alpha\| X \|_{\psi_\alpha}} \right) \le  \left(\E\phi_\alpha\left(\frac{\left|\E\left(X|\ccF \right)\right|}{\| X \|_{\psi_\alpha}}\right)\right)^{\frac{1}{c_\alpha^\alpha}} \le 2,
          \end{equation*}
          which completes the proof.
       \end{proof}
    \end{lema}
    Now we give two concentration inequalities which are valid for random variables with finite Orlicz norm. The first one is an easy consequence of the Markov inequality, therefore we omit the proof.
    \begin{lema}\label{OCI}
       For any random variable $X$ with $\|X\|_{\psi_\alpha} < \infty$ and $t> 0$ we have
       $$\P\left(|X| \ge t\right) \le 2 \exp\left(-\frac{t^\alpha}{\|X\|_{\psi_\alpha}^\alpha}\right).$$
    \end{lema}

    \begin{lema}[Tail inequality for conditional mean value] \label{TICMVL}
       Let $0 < \alpha \le 1$. Assume that a random variable $X$ satisfies $\| X \|_{\psi_\alpha} < \infty$. Moreover, let $\ccF$ be some sigma field. Then for any $t \ge \left(\frac{2}{\alpha}\right)^{1/\alpha}\| X \|_{\psi_\alpha}$ we have
       $$\P\left(\left|\E(X|\ccF)\right| > t\right) \le 6 \exp\left(-\frac{t^\alpha}{2\|X\|_{\psi_\alpha}^\alpha} \right).$$
       \begin{proof}
          Fix $c > \|X\|_{\psi_\alpha}$ and $t \ge \left(\frac{2}{\alpha}\right)^{1/\alpha}c$. Then in particular we have $\alpha \left(\frac{t}{c}\right)^\alpha \ge 2$. Using the Markov and Jensen inequalities along with $  \Gamma(x)\le  x^x / e^{x - 1} \textnormal{ (\cite{LC}, Thm. 1)}$ and  Lemma \ref{moment_estimation}  with  $Y = (|X|/c)^\alpha$,  $\beta = t^\alpha  /c^\alpha $ we get
          \begin{equation*}
             \begin{split}
                &\P\left(\left|\E(X|\ccF)\right| > t\right)  \le  \P\left(\left|\E(X|\ccF)\right|^{\alpha \frac{t^\alpha}{c^\alpha}} > t^{\alpha \frac{t^\alpha}{c^\alpha}}\right) \le  t^{-\alpha \left(\frac{t}{c}\right)^\alpha} \E\left|\E(X|\ccF)\right|^{\alpha \left(\frac{t}{c}\right)^\alpha}\le  t^{-\alpha \left(\frac{t}{c}\right)^\alpha} \E\left|X\right|^{\alpha \left(\frac{t}{c}\right)^\alpha} \\
                & = (t/c)^{-\alpha \left(\frac{t}{c}\right)^\alpha} \E\left|X/c\right|^{\alpha \left(\frac{t}{c}\right)^\alpha} \le 2e\left(t/c\right)^\alpha \exp\left(-\left(t/c\right)^\alpha\right)\le 2 e\exp\left(-(1/2)\left(t/c\right)^\alpha\right),
             \end{split}
          \end{equation*}
          where  in the last inequality we used the estimate $xe^{-x} \le e^{-\frac{x}{2}}$ which is  valid for all $x \in \R$. Now it is enough to take limit $c\rightarrow \| X \|_{\psi_\alpha}$ and notice that $2e \le 6$.
       \end{proof}
    \end{lema}

\end{document}